\documentclass{m2an-clean}
\pdfoutput=1
\usepackage{mathtools}     
\usepackage{tikzgraphicx}
\usepackage[f]{esvect}     
\usepackage{graphicx}
\graphicspath{{figures/}}

\definecolor{myred}{rgb}{0.8,0,0.1}
\definecolor{nightblue}{rgb}{0.06,0.02,0.42}
\definecolor{mygreen}{rgb}{0.0,0.54,0.11} 
\definecolor{myorange}{rgb}{0.8,0.33,0.0} 

\newcommand\psimax{\psi_\mathrm{max}}
\newcommand\phimax{\phi_\mathrm{max}}

\newcommand{\arc}{\sigma}
\newcommand{\Ucal}{\mathcal{U}}

\newcommand{\Tcal}{\mathcal{T}}

\newcommand\ysol{\bar{y}}
\newcommand\xsol{\bar{x}}
\newcommand\usol{\bar{u}}
\newcommand\zsol{\bar{z}}
\newcommand\psol{\bar{p}}
\newcommand\tfsol{\bar{t}_f}
\newcommand\tsol{\bar{t}\hspace{0.1em}}
\newcommand\etasol{\bar{\eta}}

\renewcommand\ie{\textit{i.e.}~}

\newcommand{\setR}{\mathbb{R}}
\newcommand{\hampath}{\texttt{HamPath}}
\newcommand{\bocop}{\texttt{Bocop}}
\newcommand{\matlab}{\texttt{Matlab}}

\newcommand{\refHyp}[1]{($\mathbf{A_#1}$)} 

\newcommand{\intervalle}[4]{\mathopen{#1}#2
                                \mathclose{}\mathpunct{},#3
                                \mathclose{#4}}
\newcommand{\intervalleff}[2]{\intervalle{[}{#1}{#2}{]}}
\newcommand{\intervalleof}[2]{\intervalle{(}{#1}{#2}{]}}
\newcommand{\intervallefo}[2]{\intervalle{[}{#1}{#2}{)}}
\newcommand{\intervalleoo}[2]{\intervalle{(}{#1}{#2}{)}}

\newcommand{\abs}[1]{\lvert#1\rvert} 
\newcommand{\norme}[1]{\lVert#1\rVert}

\newcommand{\fracpartial}[2]{
    \frac{\partial {#1}}{\partial {#2}}
}

\newcommand{\diff}{\mathop{}\mathopen{}\mathrm{d}}
\newcommand{\veps}{\varepsilon}

\newcommand{\ps}[2]{\left\langle #1\mathpunct{},#2\right\rangle}

\newcommand{\enstq}[2]{\left\{#1\mathrel{}\middle|\mathrel{}#2\right\}} 
\newcommand{\prodscal}[2]{\left\langle#1,#2\right\rangle} 

\DeclareMathOperator{\rank}{rank}

\newcommand{\medhrule}{\hline}
\newcommand{\bighrule}{\hline}

\makeatletter
\newcommand{\reqnomode}{\tagsleft@false\let\veqno\@@eqno}
\newcommand{\leqnomode}{\tagsleft@true\let\veqno\@@leqno}
\makeatother

\newcommand{\nbSet}[1]{\mathbb{#1}}
\newcommand{\setPositive}{\text{\bf{\tiny+}}}
\newcommand{\setNegative}{\mathbb{\tiny-}}
\newcommand{\setStar}{\text{*}}

\newcommand{\setDeco}[2]{
    \IfEqCase{#2}{
        {s}{\nbSet{#1}^{\setStar}}
        {n}{\nbSet{#1}^{\phantom{\setStar}}_{\setNegative}}
        {p}{\nbSet{#1}^{\phantom{\setStar}}_{\setPositive}}
        {sn}{\nbSet{#1}^{\setStar}_{\setNegative}}
        {sp}{\nbSet{#1}^{\setStar}_{\setPositive}}
    }
}

\newcommand{\R}{ \ensuremath{\nbSet{R}} }

\newcommand\boundaryxpos[1]{
\begin{tikzpicture}
\begin{scope}[scale = #1]
\node (A1)  at (-2,0) {};
\node (A2)  at (-1.5,0) {};
\node (A3)  at (-1.0,0) {};
\node (A4)  at (-0.5,0) {};
\node (A5)  at (0.5,0) {};
\node (A6)  at (1.0,0) {};
\node (A7)  at (1.5,0) {};
\node (A8)  at (2.0,0) {};
\node (A0)  at (0,0) {};
\node (B0)  at (00.2,0.2) {};
\node (B1)  at (-1.8,0.2) {};
\node (B2)  at (-1.3,0.2) {};
\node (B3)  at (-0.8,0.2) {};
\node (B4)  at (-0.3,0.2) {};
\node (B5)  at (0.7,0.2) {};
\node (B6)  at (1.2,0.2) {};
\node (B7)  at (1.7,0.2) {};
\node (B8)  at (2.2,0.2) {};

\node (cent1) at ($(A0) - 0.2*(0.0,1.0)$){};
\node (x0) at (-1.25,0) {};
\node (x0_) at (-1.25,-0.2) {$q_0$};
\node (x1) at (1.25,0)  {};
\node (x1_) at (1.25,-0.2)  {$q$};
\node(q)   at (-0.25, -1.1) {};

\node (x0q) at ($0.5*(q)+0.5*(x0)$) {};
\node (x1q) at ($0.5*(q)+0.5*(x1)$) {};
\node (x0x1) at ($0.45*(x1)+0.45*(x0)$) {};

\node (Sb) at (-2.0,-0.75){};
\node (Se) at ($(Sb)+2.0*(0.75,0.75)$){};
\node (SbSe) at ($(Sb)+0.8*(0.75,0.75)$){};

\draw (A1.center) -- (A8.center);
\draw (A0.center) -- (B0.center);
\draw (A1.center) -- (B1.center);
\draw (A2.center) -- (B2.center);
\draw (A3.center) -- (B3.center);
\draw (A4.center) -- (B4.center);
\draw (A5.center) -- (B5.center);
\draw (A6.center) -- (B6.center);
\draw (A7.center) -- (B7.center);
\draw (A8.center) -- (B8.center);

\draw[nightblue] (Sb)--(Se) node [pos = 0.1,right] {\hspace{0.2em}$S$};
\draw[>=latex,nightblue] (Sb)--(SbSe) ;
\draw[mygreen] (x0.center) -- (q.center) node[pos = 0.5, left] {$\arc_-$};
\draw[->, >=latex,mygreen] (x0.center) -- (x0q.center);
\draw[myred] (q.center) -- (x1.center) node[pos = 0.5, right] {~$\arc_+$};
\draw[->, >=latex,myred] (q.center) -- (x1q.center);
\draw[myorange] (x0.center) -- (x1.center); 
\draw[->, >=latex,myorange] (x0.center) -- (x0x1.center) node[pos=1.0,below] {$\arc_c$};

\node (x0) at (-1.25,0) {$\bullet$};
\node (x1) at (1.25,0)  {$\bullet$};
\node (a)  at ($(q) - (0.0,0.2)$) {$b>0$}; 
\end{scope}
\end{tikzpicture}}

\begin{document}
\title{Singular versus boundary arcs for aircraft trajectory optimization in climbing phase}
\author{Olivier Cots}\address{Toulouse Univ., INP-ENSEEIHT-IRIT, UMR CNRS 5505, 2 rue Camichel, 31071 Toulouse, France; 
    ~\\\email{olivier.cots@irit.fr\ \&\ joseph.gergaud@irit.fr\ \&\ boris.wembe@irit.fr}}
\author{Joseph Gergaud}\sameaddress{1}
\author{Damien Goubinat}\address{Thales Avionics SA, 105 av du General Eisenhower, B.P. 1147, 31047 Toulouse Cedex, France; 
    ~\\\email{damien.goubinat@fr.thalesgroup.com}}
\author{Boris Wembe}\sameaddress{1}
\date{2022}
\begin{abstract} In this article, we are interested in optimal aircraft trajectories in climbing phase. 
    We consider the cost index criterion which is a convex combination of the time-to-climb and the fuel consumption.
    We assume that the thrust is constant and we control the {flight path angle} of the aircraft. This optimization problem is 
    modeled as a Mayer optimal control problem with a single-input affine dynamics in the control and with two pure 
    state constraints, limiting the Calibrated AirSpeed (CAS) and the Mach speed. The candidates as minimizers are 
    selected among a set of extremals given by the maximum principle. We first analyze the minimum time-to-climb 
    problem with respect to the bounds of the state constraints, combining small time analysis, 
    indirect multiple shooting and homotopy methods with monitoring. This investigation emphasizes two strategies:
    the common CAS/Mach procedure in aeronautics and the classical Bang-Singular-Bang policy in control theory.
    We then compare these two procedures for the cost index criterion. \end{abstract}
\begin{resume} Dans cet article, on considère le problème d'optimisation de trajectoires d'avions en phase de montée.
Le critère à minimiser est un compromis entre la durée de la montée et la consommation de carburant. 
On suppose la poussée des moteurs constante et maximale et on contrôle la pente air. Ce problème d'optimisation 
peut s'écrire sous la forme d'un problème de contrôle optimal de Mayer dont la dynamique est affine en le contrôle.
Le problème possède deux contraintes d'états, la première sur la vitesse air calibrée (CAS) et la seconde sur
la vitesse de Mach. Le principe du maximum de Pontryagin nous donne un ensemble de conditions nécessaires d'optimalité
que l'on utilise à la fois pour analyser les solutions et pour définir les méthodes numériques de tir multiple indirect
et d'homotopie. Une analyse fine à temps court du problème de temps minimal révèle l'importance de deux types de
stratégies : la procédure classique en aéronautique appelée procédure CAS/Mach et la politique Bang-Singulière-Bang, 
classique en théorie du contrôle optimal. On compare numériquement ces deux stratégies pour le critère combinant
durée de montée et consommation de carburant.\end{resume}
\maketitle


\section{Introduction}

\label{sec:intro}
The climbing phase is a normalized phase where most of the civil aircrafts follow the CAS/Mach procedure. This procedure splits the climbing phase in two parts, the first one is an arc at constant \emph{Calibrated AirSpeed} (CAS) and the second one is an arc at constant \emph{Mach} number. The CAS constant is chosen among the ones which maximize the climbing slope, the vertical speed or the climbing rate \cite{hull:07,verriere:97}. The Mach constant depends most of the time on the cruise Mach number but this value is also closed to the one which maximizes the vertical speed. The climbing path resulting from a CAS/Mach couple may be seen as an approximation of an optimal flight, minimizing a cost made by a convex combination of the time-to-climb and the fuel consumption, that we call the \emph{cost index}. Even though historically a lot of studies was done to optimize the climbing flight path \cite{erzberger:75,miele_2:52,miele_1:52}, the CAS/Mach procedure is still used mainly because of its simplicity. The CAS and the Mach numbers are fairly easy quantities to compute since they are defined thanks to the differential pressure which is directly available from the {Pitot} tubes. 

In this article, we study first the minimum time-to-climb trajectory submitted to state constraints which are related to the CAS and the Mach number. The dynamics of this climbing phase is depicted by the following four-dimensional dynamical system:
\begin{align}
    \label{eq:systemH} \frac{\diff h}{\diff t}      &= v\, \sin \gamma                                                                  \\
    \label{eq:systemV} \frac{\diff v}{\diff t}      &= \frac{T(h)}{m} - \frac{1}{2} \frac{\rho(h) S v^2}{m} C_d(C_l) - g_0\, \sin \gamma  \\
    \label{eq:systemM} \frac{\diff m}{\diff t}      &= - C_s(v)\, T(h)                                                                  \\
    \label{eq:systemG} \frac{\diff \gamma}{\diff t} &= \frac{1}{2} \frac{\rho(h) S v}{m} \, C_l - \frac{g_0}{v} \cos \gamma,
\end{align}
where the state variable is composed of the altitude, the true air speed, the mass and the {flight path angle} of the aircraft. We refer to \cite{theseDamien} for more details about the dynamics. The altitude $h$ is given in meter (m), the true air speed $v$ in meter per second (m.s$^{-1}$), the mass $m$ in kilogram (kg) and the {flight path angle} $\gamma$ in radian (rad). In this model, the lift coefficient $C_l$ may be considered as the control variable. The BADA \cite{poles:bada09} model is chosen to represent the aircraft performance model and provides the following expressions:
\begin{equation*}
    C_s(v)   \coloneqq  C_{s,1} \left(1 + \frac{v}{C_{s,2}} \right),          \quad
    T(h)     \coloneqq C_{T,1} \left(1 - \frac{h}{C_{T,2}} + h^2 C_{T,3}\right), \quad
    C_d(C_l) \coloneqq C_{d,1} + C_{d,2}\, C_l^2, 
\end{equation*}
where the constants $C_{T,i}$, $C_{s,i}$ and $C_{d,i}$ depend on the flight phase and on the aircraft. The {International Standard Atmospheric} (ISA) model is used to represent the atmospheric model from the sea level till the end of the tropopause which here is considered at $11\,000$ meters. This model provides the evolution of the pressure $P$, the temperature $\Theta$ and the air density $\rho$ with respect to the altitude through the following expressions:
\begin{equation*}
        P(h)        \coloneqq P_0\left(\frac{\Theta(h)}{\Theta_0}\right)^{\frac{g_0}{\beta R}}, \quad
        \Theta(h)   \coloneqq \Theta_0 - \beta h, \quad
        \rho(h)     \coloneqq \frac{P(h)}{R\Theta(h)}.
\end{equation*} 
The remaining data are positive constants: $g_0$ is the gravitational constant at the sea level, $S$ the wing area, $R$ the specific constant of air, $\beta$ the thermal gradient and $P_0$, $\Theta_0$ represent the pressure and the temperature at the sea level. See Table~\ref{tab:constants} for the chosen values of the constant parameters for the numerical experiments.
\begin{table}[t!]
    \centering
    \caption{Medium-haul aircraft parameters.}
    \begin{tabular}{lll|lll}
            \medhrule
            Data       & Value     & Unit                & Data   & Value        & Unit                  \\
            \bighrule
            S          & 122.6             & m$^2$       & $C_{s,1}$      & 1.055e$^{-5}$  & kg.s$^{-1}$.N$^{-1}$ \\
            g          &   9.81            & m.s$^{-2}$  & $C_{s,2}$      & 441.54           & m.s$^{-1}$           \\
            $C_{T,1}$  & 141040            & N           & R              & 287.058          & J.kg$^{-1}$.K$^{-1}$  \\
            $C_{T,2}$  & 14909.9           & m           & $\Theta_0$     & 288.15           & K           \\
            $C_{T,3}$  & 6.997e$^{-10}$  & m$^{-2}$      & $\beta$        & 0.0065           & K.m$^{-1}$  \\
            $C_{d,1}$  & 0.0242            &             & $P_0$          & 101325           & Pa          \\
            $C_{d,2}$  & 0.0469            &             &               &                  &             \\
            \medhrule
    \end{tabular}
    \label{tab:constants}
\end{table}

Taking into account the {flight path angle} $\gamma$ in the dynamics introduces numerical instabilities \cite{ocam:2017, ifac:2017} which are known as \emph{singular perturbations} \cite{ardema:1977,survey_sing:2001,moissev:1985} and which come from large time constant differences between the state variables. A comparison of the time constants \cite{ardema:1977, theseDamien} shows that the dynamics \eqref{eq:systemH}--\eqref{eq:systemG} contains slow (the mass $m$) and fast (the {flight path angle} $\gamma$) variables. The altitude $h$ and the true air speed $v$ are fast compare to the mass but slow compare to the {flight path angle}. In this article, we consider $h$ and $v$
as slow variables. Roughly speaking, the singular perturbation phenomenon arise when a system of differential equations contains at least one small parameter multiplying the derivative of one or more state variables. To emphasize the presence of a singular perturbation in the {flight path angle} dynamics, we introduce a parameter $\veps > 0$ such that eq.~\eqref{eq:systemG} is replaced by:
\begin{equation}
    \tag{\ref{eq:systemG}'}
    \veps\, \frac{\diff \gamma}{\diff t} = \frac{1}{2} \frac{\rho(h) S v}{m} \, C_l - \frac{g_0}{v} \cos \gamma.
    \label{eq:gammaEps}
\end{equation}

Let $(P_{t_f}^{\veps})$ denote the minimum time-to-climb problem with the additional artificial parameter $\veps > 0$. From the control theory, for a fixed value of $\veps > 0$, the candidates as minimizers are selected among a set of \emph{BC-extremals}, solution of a Hamiltonian system given by the \emph{Pontryagin Maximum Principle}
(PMP) \cite{pontryaguine:1962}. The application of the PMP leads to define a Boundary Value Problem denoted (BVP$^\veps$), in terms of state and adjoint variables, which can be solved using indirect multiple shooting methods \cite{bulirsch:2002}. It is well known that multiple shooting increases numerical stability and a good alternative would be to use direct multiple shooting \cite{BOCK19841603} to solve the optimal control problem. Difficulties in solving the singularly perturbed boundary value problem
(BVP$^\veps$) arise because the solution exhibits narrow regions of very fast variation. In this case, it may be difficult to determine a mesh for the nodes of the multiple shooting method that will give an accurate numerical solution and to give an initial guess that will lead to the convergence of the underlying Newton-like algorithm. In this context, an alternative approach would be to solve the boundary value problem using automatic mesh refinement \cite{cash:2006,Cash:2001}, or to use homotopy techniques \cite{ifac:2017, CGW:2020}.

Another point of view to deal with singular perturbations is to approximate the solution using asymptotic expansions. An overview of methods to deal with singular perturbations may be found in \cite{ardema:1977,survey_sing:2001,moissev:1985}.
We consider in this article only the zero-order term of the asymptotic expansion.
To uniformly approximate the solution of a singularly perturbed boundary value problem, one has to compute at least two different approximations (the so-called inner and outer solutions) which are accurate only for part of the range of the time variable, and usually they are valid on different length-scales. The Method of Matched Asymptotic Expansion consists in computing the inner and outer solutions and then in matching them to get a uniform approximation which satisfies the boundary conditions. In this article, we consider only the zero-order approximation of the outer solution since it has been shown for the minimum time-to-climb problem with no state constraints that it is a good approximation of the solution \cite{ifac:2017}. Note that in the aircraft dynamics the zero-order reduction is equivalent to the quasi-steady approximation of the flight which is commonly performed \cite{espin:2014,hull:07,lemerrer:2012,nguyen:2006,verriere:97}.

The reduced-order optimal control problem denoted $(P_{t_f})$ is obtained putting $\veps$ to $0$ in eq.~\eqref{eq:gammaEps} and considering $\gamma$ as the new control variable. The problem $(P_{t_f})$ has one state variable less than $(P_{t_f}^{\veps})$, $\veps > 0$, and can be tackled with the tools from geometric optimal control theory.
We refer to \cite{ardema:1976,ardema:1977,theseDamien} for details about the reduction process. In particular, one important result is that the reduced-order boundary value problem obtained from (BVP$^\veps$) putting $\veps$ to $0$ is equivalent to the boundary value problem obtained after applying the maximum principle to the reduced-order dynamical system obtained from \eqref{eq:systemH}--\eqref{eq:gammaEps} putting $\veps$ to $0$, when we consider at the end $\gamma$ as the new control variable. Besides, the reduced-order dynamical system with $\gamma$ as the control variable is affine with respect to the control if we consider the small-angle approximation, that is if we replace $\cos\gamma$ by $1$ and $\sin\gamma$ by $\gamma$ assuming $\gamma$ is small. In this case, the associated pseudo-Hamiltonian is of the form $H(x,p,u) = H_0(x,p) + u\, H_1(x,p)$ where $x \coloneqq (h,v,m)$ is the reduced state, $p$ is the adjoint variable, $u$ the control (\ie the {flight path angle} $\gamma$),
and where $H_0$ and $H_1$ are two Hamiltonian lifts. In addition, to complete the definition of the optimal control problem, the control has to satisfy a constraint of the form $u \in \intervalleff{u_\mathrm{min}}{u_\mathrm{max}}$ and the state has to satisfy two constraints denoted $\phi(x) \le \phimax$ and $\psi(x) \le \psimax$, where $\phi$ represents the CAS, $\psi$ represents the Mach, and where $\phimax$ and $\psimax$ may be seen as given speed limitations.

\medskip
\paragraph*{{\textbf{Main contributions and results}}}%
The minimum time-to-climb problem with no state constraints, or from another point of view, with state constraints but with $\phimax$ and $\psimax$ big enough, is analyzed in \cite{eucass:2017,ifac:2017}. In \cite{eucass:2017}, the influence of the initial mass $m_0$ and the final true air speed $v_f$ is studied and it is shown in particular that the trajectories are of the form $\arc_\pm\arc_s\arc_\pm$, where
$\arc_-$ represents a \emph{bang} arc where $u(t) = u_\mathrm{min}$ along the arc, where $\arc_+$ represents a bang arc with $u(t) = u_\mathrm{max}$ and $\arc_s$ represents a \emph{singular} arc where $u(t) \in \intervalleoo{u_\mathrm{min}}{u_\mathrm{max}}$. The first result of this article is that these trajectories from the state unconstrained case are compliant with the so-called \emph{Maximum Operating Speed} (VMO) and \emph{Maximum Operating Mach} (MMO), that is $\phi(x(t)) < \mathrm{VMO}$ and $\psi(x(t)) < \mathrm{MMO}$
along the trajectories. Thus, with the constraints VMO and MMO, the trajectories are of the form Bang-Singular-Bang. Another contribution of this article is to analyze the influence of the bounds $\phimax$ and $\psimax$ on the structure of the trajectories for the minimum time-to-climb problem. We propose a methodology introduced in \cite{cots:2017}, based on differential homotopy \cite{hampath:2012} and geometry techniques \cite{bonnard:2006} to classify the different structures with respect to $\phimax$ and $\psimax$. Remarkably, the CAS/MACH procedure appears in our classification. In a final part, we present a comparison between trajectories that follow the simple CAS/Mach procedure and minimum cost index solutions (with VMO and MMO bounds), this criterion being commonly used in operations.
{Our numerical analysis shows that the CAS/MACH trajectories are close to the optimal solutions up to 0.3\%, on average,
regarding to the cost index criterion (see Table~\ref{tab:nlp_output} for details). These optimal solutions being of 
the form Bang-Singular-Bang.}

\medskip
The paper is organized as follows. The physical model with the statements of the optimal control problems are given in Section~\ref{sec:definitions_ocp} while the necessary conditions of optimality for both the minimum time-to-climb and the fuel consumption problems are given in Sections.~\ref{sec:boundary_def}-\ref{sec:boundary_param}. 
A small time analysis which describes the behavior of \emph{hyperbolic} trajectories with state constraints of \emph{order one} is given in Sections~\ref{subsec:classification} and \ref{subsec:small_time}.
We illustrate in Section~\ref{subsec:cas_mach_make_sense}, how this analysis reveals the CAS/Mach procedure. 
This analysis is then applied to the minimum time-to-climb problem in Section~\ref{sec:minimum_time_application} combined with numerical methods (indirect multiple shooting and homotopy with monitoring) to classify the optimal structures with respect to the bounds $\phimax$ and $\psimax$.
A comparison between minimal cost index solutions (with singular arcs) and the CAS/Mach procedure is then proposed in Section~\ref{sec:cost_index_and_cas_mach}. Finally, Section~\ref{sec:conclusion} concludes the article.

\section{Preliminary general results and theoretical background}
\label{sec:general}

\subsection{Statement of the optimal control problems}
\label{sec:definitions_ocp}

We restrict the dynamics to the vertical motion of the aircraft. A complete description of the motion can be found in \cite{theseDamien}. The aircraft is subjected to four forces, the Drag $\vv{D}$, the Lift $\vv{L}$, the Thrust $\vv{T}$ and its own weight $\vv{P}$. A nonlinear point mass representation is used and we consider that all the forces apply on the center of gravity of the aircraft. We assume that the thrust is collinear to the velocity vector $\vv{V}$, which means that the angle of attack is neglected here, and that the aircraft evolves in an horizontal constant wind field. The application of the first dynamics principle assuming that the earth is Galilean, provides the four eqs.~\eqref{eq:systemH}--\eqref{eq:systemG}.

During the flight, the lift coefficient $C_l$ depends on the variation of the angle of attack. Since this quantity is not taken into account here, it is quite natural to consider the lift coefficient $C_l$ as the control variable. However, according to the work presented in \cite{ardema:1977,survey_sing:2001} this four-dimensional
dynamics contains slow ($m$) and fast ($\gamma$) variables. The time scale separation between the slow and fast variables is handled by a singular perturbation analysis which consists here in computing a zero-order approximation of the solution. The reduction process has two main features. In a first step, we substitute eq.~\eqref{eq:systemG} by the quasi-steady approximation:
\begin{equation*}
0 = \frac{1}{2}\frac{\rho S v}{m} C_l -\frac{g_0}{v}\cos\gamma
%
~\Rightarrow~
%
C_l = \frac{2 m g_0}{\rho S v^2}\cos\gamma.
\end{equation*}
In a second time, the lift coefficient is replaced by the previous expression in eqs.~\eqref{eq:systemH}--\eqref{eq:systemM} and the {flight path angle} $\gamma$ is taken as the new control variable. We also consider that the {flight path angle} remains small ($\gamma$ varies from $u_{\min}\coloneqq 0$ to $u_{\max}\coloneqq 0.262$ rad) during the climbing so we set $\cos\gamma \approx 1$ and $\sin\gamma \approx\gamma$. These considerations lead to the new affine control system:
\begin{equation}
    \label{eq:sys_control}
\dot{x}(t) = F_0(x(t)) + u(t)\, F_1(x(t)),
\end{equation}
where $x\coloneqq (h,v,m)$, $u\coloneqq \gamma$,
\begin{equation*}
    F_0(x) \coloneqq 
    \left(\begin{aligned}
    &0 \\
    \frac{T(h)}{m} -\frac{1}{2}\frac{\rho(h)Sv^2}{m}&C_{d,1} - 2\frac{mg_0^2}{\rho(h)Sv^2}C_{d,2}\\
    - C_s&(v)T(h)
    \end{aligned}\right)
    %
    \quad \text{and} \quad
    %
    F_1(x) \coloneqq \left(\begin{aligned}
    &v\\
    -&g_0\\
    &0
    \end{aligned}\right).
\end{equation*}
According to \cite{ifac:2017}, this reduced dynamics is a sufficient approximation of the initial dynamics considering the state unconstrained minimum time-to-climb problem.

The true air speed of an aircraft is quite difficult to measure during the flight, that is why historically, the concept of CAS was introduced. This speed can be computed using the differential ratio of pressure and can be expressed, through Bernoulli's equations, as a function of the true air speed. The Mach speed, for its part, is defined by the ratio between the true air speed and the speed of the sound. Here, the CAS is denoted by $\phi$ and the Mach by $\psi$ and the corresponding expressions are given by:
\begin{align*}
\phi(x) & \coloneqq \sqrt{2\frac{R\,\Theta_0}{\kappa}
    \left(
    \left(
        \frac{P(h)}{P_0}
    \left(
        \frac{\kappa\, v^2}{2R\,\Theta(h)}+1
    \right)^{\kappa^{-1}}
        +1
    \right)^{\kappa}-1
    \right)}, 
\\
\psi(x) &\coloneqq \frac{v}{\sqrt{\gamma_\mathrm{air}R\,\Theta(h)}},
\end{align*}
where $\gamma_\mathrm{air}\coloneqq 1.4$ J.K$^{-1}$ is the heat capacity of the air and $\kappa\coloneqq \gamma_\mathrm{air}/(1-\gamma_\mathrm{air})$. From these expressions, we define the state constraints $c_1(x) \le 0$ and $c_2(x) \le 0$ with:
\begin{align}
c_1(x) &\coloneqq \phi(x) - \phi_{\max},\label{eq:c1}\\
c_2(x) &\coloneqq \psi(x) - \psi_{\max},\label{eq:c2}
\end{align}
where $\phi_{\max}$ (resp. $\psi_{\max}$) can be fixed to VMO (resp. MMO) or to smaller specified values.

In this article we are interested in the minimization of different criteria: the final time $t_f$, the fuel consumption $\Delta_m \coloneqq m_0-m_f$, where $m_0$ is the given initial mass and $m_f$ the final mass, and a convex combination of the final time and the fuel consumption: $g_\alpha(t_f, m_f) \coloneqq \alpha\, t_f + (1-\alpha) (m_0-m_f).$ The criterion $g_\alpha$, for $\alpha \in \intervalleff{0}{1}$, is called the cost index. The initial state will be fixed to the realistic values $x_0 \coloneqq (h_0, v_0, m_0) \coloneqq (3480,128.6,69000)$. We define the set $\Ucal_t$ of \emph{admissible controls}, that is the subset of 
$
\enstq{u}{u\colon \intervalleff{0}{t} \to \intervalleff{u_{\min}}{u_{\max}} \text{ measurable}}
$
such that the corresponding trajectory solution of eq.~\eqref{eq:sys_control} with $x(0) = x_0$ is well defined over $\intervalleff{0}{t}$. The \emph{minimum cost index problem} can be summarized by:
\leqnomode
\begin{equation}
\tag{$P_\alpha$}
\left\{
\begin{aligned}
& \underset{(t_f, u) \in D}{\min}~ g_\alpha(t_f, m(t_f)), \\[0.5em]
& \dot{x}(t) = F_0(x(t)) + u(t)\, F_1(x(t)), \quad u(t) \in U, \quad
t \in \intervalleff{0}{t_f} \text{ a.e.}, \quad x(0)=x_0,\\[0.5em]
& c_1(x(t)) \le 0, \quad t \in \intervalleff{0}{t_f},\\
& c_2(x(t)) \le 0, \quad t \in \intervalleff{0}{t_f},\\
&b(x(t_f)) = 0,
\end{aligned}
\right.
\label{eq:ocp_cost_index}
\end{equation}
\reqnomode
where $U \coloneqq \intervalleff{u_{\min}}{u_{\max}}$, $D \coloneqq \enstq{(t_f, u)}{t_f \ge 0,~ u \in \Ucal_{t_f}}$ and $b(x) \coloneqq (h-h_f, v-v_f)$
%
%
with $h_f \coloneqq 9144$ and $v_f \coloneqq 191$. Thus, the final altitude and true air speed are fixed while the final mass is free. The \emph{minimum time-to-climb problem}, denoted ($P_{t_f}$), is defined as (\ref{eq:ocp_cost_index}$|_{\alpha=1}$)
while the \emph{minimum fuel consumption problem}, denoted ($P_{\Delta_m}$), is defined as (\ref{eq:ocp_cost_index}$|_{\alpha=0}$).

\medskip
\paragraph*{\textbf{Notation.}} %
In this paper, we follow the presentation of \cite{bonnard:2006}, which exhibits the role of the Lie and Poisson brackets that we define below.
Let $F_0$, $F_1$ be two smooth vector fields and $c$ a scalar smooth function on $\R^n$. The \emph{Lie derivative} of $c$ along $F_0$ denoted $F_0 \cdot c$ is simply the directional derivative of $c$ at $x$ along $F_0(x)$, given by $(F_0 \cdot c)(x) \coloneqq c'(x)\, F_0(x)$. The \emph{Lie bracket} between $F_0$ and $F_1$ is given by $[F_0,F_1] \coloneqq  F_0 \cdot F_1 - F_1 \cdot F_0$, with $(F_0 \cdot F_1)(x) \coloneqq F_1'(x) \, F_0(x)$. Denoting $p$ the adjoint variable and denoting $H_0(x,p) \coloneqq \prodscal{ p }{ F_0(x) }$ and $H_1(x,p) \coloneqq \prodscal{ p }{ F_1(x) }$ the Hamiltonian lifts of $F_0$ and $F_1$, then the \emph{Poisson bracket} of $H_0$ and $H_1$ is given by $\{H_0,H_1\}  \coloneqq  \vv{H_0}\cdot H_1,$ where $\vv{H_0} \coloneqq (\fracpartial{H_0}{x}, -\fracpartial{H_0}{p} )$ is the \emph{Hamiltonian system} associated to $H_0$. We also use the notation $H_{01}$ (resp. $F_{01}$) to write the bracket $\{H_0,H_1\}$ (resp. $[F_0,F_1]$) and so forth. Besides, since $H_0$ and $H_1$ are two Hamiltonian lifts, we have $\{H_0,H_1\}= \prodscal{ p }{ [F_0,F_1] }$.

\begin{rmrk}
See \cite[Chapter 2]{theseDamien} for the expressions of $F_{01}$, $F_{001}$ and $F_{101}$ 
associated to eq.~\eqref{eq:sys_control}.
\end{rmrk}



\subsection{Boundary arcs and assumptions}
\label{sec:boundary_def}

We consider in the following of Section~\ref{sec:general} the problem \eqref{eq:ocp_cost_index} but with only one state constraint, since, according to the numerical results from the next sections, the constraints $c_1$ and $c_2$ are active at the same time only at isolated times. In other words, we never have $c_1(x(t))=c_2(x(t))=0$ along a time interval of non-empty interior. Due to these considerations, we present the necessary conditions of optimality considering only one scalar state constraint denoted $c$. The optimal control problem may be written in the form:
\begin{equation*}
    \begin{aligned}
        \underset{(t_f, u) \in D}{\min}& \quad g_\alpha(t_f, m(t_f, x_0, u)) \\
        \text{subject to } & \\
        & b(x(t_f, x_0, u)) = 0, \\
        & c(x(t, x_0, u)) \le 0,~ \forall\, t \in \intervalleff{0}{t_f},
    \end{aligned}
\end{equation*}
where $x(\cdot,x_0,u)$ is the solution of eq.~\eqref{eq:sys_control} with the initial condition $x(0,x_0,u) = x_0$.

We call a \emph{boundary arc}, labeled $\arc_c$, an arc defined on an interval $I\coloneqq\intervalleff{a}{b}$ (not reduced to a singleton), such that $c(\arc_c(t)) = 0$, for every $t\in I$. The times $a$ and $b$ are called the \emph{entry-} and \emph{exit-time} of the boundary arc; $a$ and $b$ are also termed \emph{junction times}. An arc $\arc$ is said to have a \emph{contact point} with the boundary at $\tsol \in \intervalleff{0}{t_f}$ if $c(\arc(\tsol))=0$ and $c(\arc(t))<0$ for $t\ne\tsol$ in a neighborhood of $\tsol$. A subarc $\arc$ on which $c(\arc(t)) < 0$ is called an \emph{interior arc}.

The \emph{generic order} of the constraint $c$ is the integer $m$ such that
$ F_1 \cdot c = F_1 \cdot (F_0 \cdot c) = \dots = F_1 \cdot (F_0^{m-2} \cdot c) = 0$
and $F_1 \cdot (F_0^{m-1} \cdot c) \ne 0$.
If the order of a boundary arc $\arc_c$ is $m$, then its associated feedback control can be generically computed by differentiating $m$ times the mapping $t \mapsto c(\arc_c(t))$ and solving with respect to $u$ the linear equation:
\[
c^{(m)} = F_0^m \cdot c + u \, F_1 \cdot (F_0^{m-1} \cdot c) = 0.
\]
The boundary feedback control denoted $u_c$ is given by
\[
u_c  \coloneqq  - \frac{F_0^m \cdot c }{ F_1 \cdot (F_0^{m-1} \cdot c)}.
\]

Let $t\mapsto\arc_c(t)$, $t\in\intervalleff{t_1}{t_2} \subset \intervalleff{0}{t_f}$,
be a boundary arc associated to $u_c(\cdot)$. We introduce the assumptions:
\begin{itemize}
\item[] \hspace{-2em} {($\mathbf{A_1}$)}  $(F_1 \cdot (F_0^{m-1} \cdot c))({\arc_c(t)}) \ne 0$ for every $t\in\intervalleff{t_1}{t_2}$, with $m$ the order of the constraint.
\item[] \hspace{-2em} {($\mathbf{A_2}$)} ${u_c(t)} \in \intervalleff{u_{\min}}{u_{\max}}$ for $t\in\intervalleff{t_1}{t_2}$, \ie the boundary control is admissible.
\item[] \hspace{-2em} {($\mathbf{A_3}$)} ${u_c(t)}\in \intervalleoo{u_{\min}}{u_{\max}}$ for $t\in\intervalleoo{t_1}{t_2}$, \ie $u_c$ is not saturating on $\intervalleoo{t_1}{t_2}$.
\end{itemize}

\begin{rmrk}
These assumptions are numerically checked \textit{a posteriori}.
\end{rmrk}

\subsection{The maximum principle}

We recall hereinafter the necessary conditions due to \cite{jacobson:1971,maurer:1977}.
Firstly, we define for $(x,p,u,\eta) \in H \colon R^n \times (R^n)^* \times \R \times \R$ the \emph{pseudo-Hamiltonian}:
\begin{equation*}
H(x,p,u,\eta)  \coloneqq  \prodscal{ p }{ F_0(x) + u\, F_1(x) } + \eta\, c(x)
\end{equation*}
where $\eta$ is the \emph{Lagrange multiplier} of the constraint $c$ and $n\coloneqq 3$ is the dimension of the state. Now, let us consider $(\tfsol, \usol(\cdot)) \in D$ an optimal solution with associated trajectory $\xsol(\cdot)$. Let us assume that the set of contact and junction times with the boundary, denoted $\Tcal$, is finite. Assume also that the optimal control is piecewise smooth and that along each boundary arc, assumptions \refHyp{1} and \refHyp{2} are satisfied.
Then, we have the following necessary optimality conditions:

\begin{enumerate}
\item There exists a function $\etasol(\cdot) \le 0$, a real number $p^0\le 0$ and a function of bounded variation $\psol(\cdot) \in BV(\intervalleff{0}{\tfsol},(\R^n)^*)$ such that for $t \in \intervalleff{0}{\tfsol} \text{ a.e.}$:
\begin{align*}
    \dot{\xsol}(t) =  \frac{\partial H}{\partial p}(\xsol(t),\psol(t),\usol(t),\etasol(t)), \quad 
    \dot{\psol}(t) = -\frac{\partial H}{\partial x}(\xsol(t),\psol(t),\usol(t),\etasol(t)).
\end{align*}
\item The \emph{maximization condition} holds for $t \in \intervalleff{0}{\tfsol} \text{ a.e.}$:
\begin{equation}
    \displaystyle H(\xsol(t),\psol(t),\usol(t),\etasol(t)) = \max_{u\in U}H(\xsol(t),\psol(t),u,\etasol(t)).
    \label{eq:maxPbConstraint}
\end{equation}
\item The \emph{boundary conditions} $b(\xsol(\tfsol))=0$ are satisfied and we have the \emph{transversality conditions}
%
    $\psol_m(\tfsol) = - p^0 (1 - \alpha)$
%
and since $t_f$ is free, if $\usol(\cdot)$ is continuous at time $\tfsol$, then,
%
 $H(\xsol(\tfsol),\psol(\tfsol),\usol(\tfsol),\etasol(\tfsol)) = - p^0 \alpha$.
%
\item The function $\etasol(\cdot)$ is continuous on the interior of the boundary arcs and we have the \emph{complementarity condition}
$
 \etasol(t)\, c(\xsol(t)) = 0$, $\forall t \in \intervalleff{0}{\tfsol}.
$
\item For any $\tau \in \Tcal$ we have
%
    $H[\tau^+] = H[\tau^-]$,
%
where $[\tau]$ stands for $(\xsol(\tau),\psol(\tau),\usol(\tau),\etasol(\tau))$, and
%
    $\psol(\tau^+)    = \psol(\tau^-) - \nu_\tau \, c'(\xsol(\tau))$,
%
where $\nu_\tau  \le 0$ is called a \emph{jump}.
\end{enumerate}

\begin{rmrk}
Either $p^0=0$ (abnormal case), or $p^0$ can be set to $-1$ by homogeneity (normal case). We restrict our study to the normal case.
\end{rmrk}

\begin{dfntn}~
\begin{itemize}
\item We call an \emph{extremal} a quadruple $(x(\cdot),p(\cdot),u(\cdot),\eta(\cdot))$ satisfying items 1, 2, 4 and 5. It is called a \emph{BC-extremal} if it satisfies also item 3.
\item Let $H_0(x,p)  \coloneqq  \prodscal{ p }{ F_0(x) }$ and $H_1(x,p)  \coloneqq  \prodscal{ p }{ F_1(x) }$ denote the {Hamiltonian lifts} of $F_0$ and $F_1$. Then, along any extremal, we call $\Phi(t)  \coloneqq  H_1(x(t),p(t))$ the \emph{switching function}.
\item It follows from \eqref{eq:maxPbConstraint} that along any extremal, we have $u(t) = u_{\min}$ if $\Phi(t) < 0$ and $u(t) = u_{\max}$ if $\Phi(t) > 0$. We say that a trajectory $x(\cdot)$ restricted to a sub-interval $I \subset \intervalleff{0}{t_f}$, not reduced to a singleton, is a \emph{bang arc} if $u(\cdot)$ is constant on $I$, taking values in $\{u_{\min},u_{\max}\}$. The trajectory is called \emph{bang-bang} if it is the concatenation of a finite number of bang arcs. 
\item We say that a trajectory $x(\cdot)$ restricted to a sub-interval $I \subset \intervalleff{0}{t_f}$, not reduced to a singleton, is a \emph{singular arc} if it is an interior arc and if the associated extremal lift satisfies $\Phi(t) = 0$, for every $t \in I$.
\item  We say that an extremal is a \emph{bang}, \emph{singular}, \emph{boundary} or \emph{interior} extremal if the associated trajectory is respectively a {bang}, {singular}, {boundary} or {interior} arc.
\end{itemize}
\end{dfntn}

\begin{dfntn}
A bang arc such that $u(\cdot) \equiv u_{\min}$ (resp. $u_{\max}$) is labeled $\arc_-$ (resp. $\arc_+$). A singular arc is labeled $\arc_s$ while
a boundary arc associated to the constraint $c$ is labeled $\arc_c$. Besides, we denote by $\arc_1\arc_2$ an arc $\arc_1$ followed by an arc $\arc_2$.
\end{dfntn}

\begin{rmrk}
Along a boundary arc, the maximization condition \eqref{eq:maxPbConstraint} with assumption \refHyp{3} imply $\Phi = 0$ on the interior of the boundary arc. Besides, the adjoint vector may be discontinuous at $\tau \in \mathcal{T}$.
\end{rmrk}

\subsection{Parameterization of the singular extremals}

Relaxing the control bounds, singular trajectories are parameterized by the constrained Hamiltonian system:
\begin{equation*}
\dot{x} = \partial_p H,\quad \dot{p} = -\partial_x H,\quad 0 = \partial_u H = H_1,
\end{equation*}
with $H_1(x,p) = \prodscal{p}{F_1(x)}$ the Hamiltonian lift of $F_1$. The constraint $H_1 = 0$ has to be differentiated at least twice along a singular extremal to compute the control. This gives:
\begin{equation*}
H_1 = H_{01} = H_{001} + u\, H_{101} = 0,
\end{equation*}
along any singular extremal. A singular extremal along which $H_{101} \ne 0$ is called of \emph{minimal order} and the corresponding singular control is given by:
\begin{equation*}
u_s(z) \coloneqq -\frac{H_{001}(z)}{H_{101}(z)},
\end{equation*}
with $z\coloneqq (x,p)$. In this case, we have the following additional necessary condition of optimality deduced from the high-order maximum principle \cite{Krener:1977}. If the singular control is not saturating along the singular extremal, then the \emph{generalized Legendre-Clebsch condition} must hold along the singular extremal, that is:
\begin{equation}
\fracpartial{}{u}\fracpartial{ {}^2 }{t^2}\fracpartial{H}{u} = H_{101} \ge 0.
\label{eq:GLC}
\end{equation}
Besides, we have the following well-known result \cite[Prop.~21]{bonnard:2003}
that we use to define the numerical methods, see Section~\ref{sec:multiple_shooting}.
\begin{prpstn}
\label{prop:singular_extremal}
Assume the open subset $\Omega \coloneqq \enstq{z}{H_{101}(z)\ne 0}$ is not empty and let us define on $\Omega$ the Hamiltonian $H_s(z) \coloneqq H_0(z) + u_s(z)\, H_1(z)$. Then, the singular extremals of minimal order are the solutions of $\dot{z}(t) = \vv{H_s}(z(t))$, starting from the set $\enstq{z}{H_1(z) = H_{01}(z) = 0}$.
\end{prpstn}

For three-dimensional systems, {the singular control may be written in feedback form only with respect to the state $x$}.
Indeed,setting 
$
D_{\xi}(x) \coloneqq \det(F_{1}(x),F_{01}(x),F_\xi(x)),
$
the singular control is given in feedback form by
\begin{equation*}
u_s(x) = -\frac{D_{001}(x)}{D_{101}(x)},
\end{equation*}
whenever $D_{101}(x)\ne 0$ and $p\ne 0$, since, along a singular extremal 
\begin{equation*}
\prodscal{p}{F_1(x)} = \prodscal{p}{F_{01}(x)} = \prodscal{p}{F_{001}(x)+u\,F_{101}(x)} = 0.
\end{equation*}
%
%

\begin{rmrk}
Assuming $D_{0}(x) \ne 0$, the generalized Legendre-Clebsch condition \eqref{eq:GLC} becomes $D_{0} D_{101} \ge 0$ (when $\alpha \ne 0$) and on the set $H_1 = H_{01} = 0$ we have $D_{101} = 0 \Rightarrow H_{101} = 0$ and $(H_{101} = 0 \text{ and } p \ne 0) \Rightarrow D_{101} = 0$.
\end{rmrk}

\subsection{Parameterization of the boundary extremals}
\label{sec:boundary_param}

We may find in \cite{bonnard_art:2003,maurer:1977} the determination of the multiplier $\eta$ and the jump $\nu_\tau$
together with the analysis of the junction conditions, which is based on the concept of order and related to the 
classification of extremals. 
{%
We give next some results only for the case $m=1$ since the constraints $c_1$ and $c_2$ are of order 1, according
to numerical experiments that we have realized as follows.

First, the differentiation of $c_1(x(t)) = \phi(x(t))-\phi_{\max}$ with respect to the time $t$ leads to 
\begin{equation*}
\dot{c}_1(x(t)) =\frac{\partial \phi}{\partial h}(x(t)) \dot{h}(t) + \frac{\partial \phi}{\partial v}(x(t)) \dot{v}(t)
%
\quad \text{and} \quad
%
\frac{\partial \dot{c}_1}{\partial u}(x) = v\frac{\partial \phi}{\partial h}(x)-g_0\frac{\partial \phi}{\partial
v}(x) = (F_1 \cdot c_1) (x).
\end{equation*}
Let us define $V^*:=140$ m.s$^{-1}$ as the minimal authorized speed in our case study\footnote{$V^*$ is closed to the maximal speed below $3480$ m: $148.5$ m.s$^{-1}$}.
Then, we check thanks to the software Sage \cite{sage:2013}, that $(F_1 \cdot c_1)(x)$ does not vanish for any $(h, v) \in \intervalleff{3400}{11000} \times \intervalleff{V^*}{VMO}$, so, for our case study, $c_1$ is a state constraint of order 1. 
The same reasoning holds for $c_2$.%
}

\begin{rmrk}
We refer to \cite[Chapter 3]{theseDamien} for the exact and quite long expressions of $F_{0}\cdot c$, $F_1 \cdot c$ and $H_{01}$ for both constraints $c_1$ and $c_2$.
\end{rmrk}

For a first-order constraint denoted $c$, assuming \refHyp{1} and \refHyp{3}, we have the following result from \cite{bonnard_art:2003}.

\begin{prpstn}
\label{prop:order1}

Let $m=1$. Then:
\begin{enumerate}
\item along the boundary, the control and the multiplier are given by 
\[
    \displaystyle u_c(x)  = -\frac{(F_{0}\cdot c) (x)}{(F_1 \cdot c) (x)} \quad \text{and} \quad
    \displaystyle \eta_c(z) = \frac{H_{01}(z)}{(F_1 \cdot c)(x)},
    \quad z = (x,p).
\]
\item if the control is discontinuous at a contact or a junction time $\tau$ between a bang arc and the boundary then the jump $\nu_\tau = 0$.
\item we have
\begin{equation*}
\nu_\tau =  \frac{\Phi(\tau^-)}{(F_{1} \cdot c) (x(\tau))} \text{ at an entry point}
%
\quad\text{and}\quad
%
\nu_\tau = -\frac{\Phi(\tau^+)}{(F_{1} \cdot c) (x(\tau))} \text{ at an exit point.}
\end{equation*}
\end{enumerate}
\end{prpstn}

Likewise the singular case, in the state constrained case, we have the following result excerpted from \cite[Prop.~4.5]{cots:2017} which is useful to define the numerical shooting method, see Section~\ref{sec:multiple_shooting}.

\begin{prpstn}
\label{prop:boundary_extremal}
Let $c(x) \le 0$ be a smooth scalar state constraint of order 1 and assume the open subset $\Omega \coloneqq \enstq{x}{ (F_1 \cdot c)(x) \ne 0}$ is not empty. On $\Omega$ we define the Hamiltonian 
\[
H_c(z) \coloneqq H_0(z) + u_c(x)\, H_1(z) + \eta_c(z)\, c(x), \quad z=(x,p),
\]
where $u_c$ and $\eta_c$ are defined in Prop.~\ref{prop:order1}. Then, the boundary extremals of order 1 contained in $H_1 = 0$ are the solutions of the system $\dot{z}(t) = \vv{H_c}(z(t))$, starting from $\enstq{z=(x,p)}{ c(x) = H_{1}(z) = 0 }$.
\end{prpstn}

\section{Small time analysis of the minimum time-to-climb problem}
\label{sec:time-to-climb}

The minimum time-to-climb problem with no state constraints is analyzed in \cite{eucass:2017,ifac:2017}. In \cite{eucass:2017}, the influence of the initial mass $m_0$ and the final true air speed $v_f$ is studied and it is shown in particular that the trajectories are of the form $\arc_\pm\arc_s\arc_\pm$, where $\arc_s$ is an {hyperbolic} singular extremal. We thus restrict the small time analysis, presented in Section~\ref{subsec:small_time}, to hyperbolic singular extremals, that we introduce in Section~\ref{subsec:classification}.



\subsection{Generic classification of extremals near the switching surface}
\label{subsec:classification}

Along a singular extremal we have $\partial_u H\equiv 0$. In our particular case of a single-input affine control system this condition becomes $\Phi(t) = H_1(z(t)) =0$. Let us define the switching surface
\[
\Sigma \coloneqq \enstq{z\coloneqq(x,p)\in \setR^n\times(\setR^n)^*}{H_1(z)=0},
\]
and the set
\[
\Sigma_s \coloneqq \enstq{z\coloneqq(x,p)\in \setR^n\times(\setR^n)^*}{H_1(z)=H_{01}(z) =0}
\]
containing all the singular extremals. Let us introduce the notation $\Phi_+$ (resp. $\Phi_-$) if the control along a bang extremal is $u_{\max}$ (resp. $u_{\min}$).
The first and second derivatives of $\Phi_\pm$ are given by:
\begin{align*}
\dot{\Phi}_\pm(t) &= H_{01}(z(t)),\\
\ddot{\Phi}_+(t)&= H_{001}(z(t)) + u_{\max}\, H_{101}(z(t)), \\
\ddot{\Phi}_-(t)&= H_{001}(z(t)) + u_{\min}\, H_{101}(z(t)).
\end{align*}
A crucial point to analyze the minimum time-to-climb problem, is then to apply the results from \cite{kupka:1987} (see also \cite{bonnard:2003}) to classify the extremal curves near the switching surface $\Sigma$. We have the following:

\medskip
\begin{enumerate}
\item \textbf{Ordinary switching time.}
It is a time $\bar{t}$ such that two bang arcs switch with $\Phi(\bar{t}) = 0$ and
$\dot{\Phi}(\bar{t}) = H_{01}(z(\bar{t})) \ne 0$. According to the maximum principle, near $\Sigma$, the extremal is of the form $\arc_-\arc_+$ if
$\dot{\Phi}(\bar{t}) > 0$ and $\arc_+\arc_-$ if $\dot{\Phi}(\bar{t}) < 0$.

\medskip
\item \textbf{Fold point.}
It is a time where a bang arc has a contact of order 2 with $\Sigma$.
We have three cases (if $\ddot{\Phi}_\pm \ne 0$) depending on $\ddot{\Phi}_\pm$ at the switching time:
\begin{itemize}
\item \emph{Hyperbolic case:} $\ddot{\Phi}_+ > 0$ and $\ddot{\Phi}_- < 0$. A connection with a singular extremal is possible at $\Sigma_s$ and locally each extremal is of the form $\arc_\pm\arc_s\arc_\pm$ (by convention each arc of the sequence can be empty).
\item \emph{Parabolic case:} $\ddot{\Phi}_+ \ddot{\Phi}_- > 0$. The singular extremal at the switching point is not admissible and every extremal curve is locally bang-bang with at most two switchings, \ie $\arc_+\arc_-\arc_+$ or $\arc_-\arc_+\arc_-$.
\item \emph{Elliptic case:} $\ddot{\Phi}_+ < 0$ and $\ddot{\Phi}_- > 0$. A connection with a singular arc is not possible and locally each extremal is bang-bang but with no uniform bound on the number of switchings.
\end{itemize}
\end{enumerate}

When dealing with three-dimensional systems we can give conditions depending only on the state $x$ to classify the three cases around fold points. Assuming $D_0(x)\neq 0$, then the family $(F_0(x),F_1(x),F_{01}(x))$ forms a basis of $\setR^3$ and there exists $(\alpha_0,\alpha_1,\alpha_{01}) \in \setR^3$ and $(\beta_0, \beta_1, \beta_{01})\in \setR^3$ such that 
\begin{align*}
F_{001}(x)+u_{\min}\,F_{101}(x) = \alpha_0 F_0(x) + \alpha_1 F_1(x) + \alpha_{01} F_{01}(x),\\
F_{001}(x)+u_{\max}\,F_{101}(x) = \beta_0 F_0(x) + \beta_1 F_1(x) + \beta_{01} F_{01}(x).
\end{align*} 
By linearity of the determinant, we have 
\begin{align*}
D_{001}(x)+u_{\min}\,D_{101}(x) = \alpha_0\, D_0(x),\\
D_{001}(x)+u_{\max}\,D_{101}(x) = \beta_0\, D_0(x),
\end{align*}
and since $D_0(x) \ne 0$ then we can compute $\alpha_0$ and $\beta_0$ along the trajectory. Besides, along any singular extremal, we have:
\begin{align*}
\ddot{\Phi}_-(t) & = \ps{p(t)}{F_{001}(x(t))+u_{\min}\,F_{101}(x(t))} = \alpha_0(t)\, H_0(z(t)), \\
\ddot{\Phi}_+(t) & = \ps{p(t)}{F_{001}(x(t))+u_{\max}\,F_{101}(x(t))} = \beta_0(t)\, H_0(z(t)),
\end{align*} 
with $H_0(z(t)) >0$ in the normal case with the convention of the maximum principle.
Denoting $\tsol$ the time when the extremal has a contact of order two with the switching surface $\Sigma$, then the singular extremal is hyperbolic if $\alpha_0(\tsol) <0$ and $\beta_0(\tsol)>0$, is elliptic if $\alpha_0(\tsol) >0$ and $\beta_0(\tsol)<0$, and parabolic if $\alpha_0(\tsol)\, \beta_0(\tsol)>0$.

\begin{rmrk}
    {%
    On \cite[Fig. 2.11, p. 59]{theseDamien}, the numerical analysis is illustrated and it shows \textit{a posteriori}
    that $\alpha_0(t) >0$ and $\beta_0(t) < 0$ along the singular arc. 
    Thus, we are in the hyperbolic case and the assumption that $D_0$ does not vanished is fulfilled.}
\end{rmrk}

\subsection{Analysis of the hyperbolic case with a one-order state constraint}
\label{subsec:small_time}

The time-optimal trajectories are concatenations of bang, singular and boundary arcs. One of the main difficulties is to determine the number and the sequence of these arcs. We can see these concatenations as a sequence of patterns (or patches), each pattern being made with only few concatenations. These patterns depend on the classification presented in Section~\ref{subsec:classification} but also on the order of the state constraint. Since we are dealing with order one state constraints and hyperbolic singular extremals, the
number of possible patterns is drastically reduced. Besides, a (long) time-optimal trajectory may be seen as a concatenation of small time-optimal trajectories since it is necessarily optimal on every sub-interval of time. Finally, the pattern we are looking for are thus given by small time-optimal trajectories that we analyze in this section.

\begin{rmrk}
    {%
    In this section, the main result is given by Prop.~\ref{prop:tps_courts} which explains the structures of the time-optimal
    trajectories. In the following, we make a brief and quite technical summary of some results excerpted
    from \cite{bonnard_art:2003,bonnard:2006}, in order to introduce the notations of Prop.~\ref{prop:tps_courts}, and which can
    be skipped in a first read.}
\end{rmrk}

\begin{rmrk}
In \cite{eucass:2017}, it is shown that the state-unconstrained time-optimal trajectories are of the form $\arc_\pm\arc_s\arc_\pm$, \ie
the optimal sequences are made with only one single pattern.
\end{rmrk}

The aircraft dynamics is a multi-scale system and since the time constant of the mass is a hundred, resp. a thousand, times greater than the time constant of the altitude, resp. the speed, we can neglect the evolution of the mass when considering short time intervals and thus, in this case, we can reduce our three dimensional dynamics into a planar dynamics. We present hereinafter some results from \cite{bonnard_art:2003,bonnard:2006} about small time analysis considering only a planar single-input control system, in the hyperbolic case and with an order one scalar state constraint $c(q(t))\le 0$, of the form:
\begin{equation*}
\dot{q}(t) = F_0(q(t)) + u(t)\, F_1(q(t)),
\end{equation*}
where $|u(t)|\le 1$ and $q(t) \coloneqq (x(t), y(t)) \in \setR^2$. Let us take $q_0$ such that $c(q_0) = 0$ and let us identify $q_0$ to the origin. Let us assume that $F_0(q_0)$ and $F_1(q_0)$ are linearly independent, and that the constraint $c$ is of order one, \textit{i.e.} $(F_1\cdot c)(q_0) \neq0$. Then replacing, if necessary, $u$ by $-u$, we can find a local diffeomorphism preserving $q_0=0$ and transforming the constrained system into the following generic model:
\[
\dot{x}(t) = 1 + y(t)\, a(q(t)),\quad \dot{y}(t) = b(q(t)) + u(t),\quad y(t) \le 0.
\]

To describe the hyperbolic case, we consider also that $\det(F_1(q_0), F_{01}(q_0)) = a(q_0) = 0$ and we assume that the set 
\[
S \coloneqq \enstq{q \in \R^2}{\det(F_1(q), F_{01}(q)) = 0}
\]
containing the singular trajectories is a simple curve (this assumption is not restrictive since our approach is local) that we approximate by a straight line in our small time model. The equations of the system become:
\[
\dot{x}(t) = 1 + y(t) \, (a\, y(t) + b\, x(t)),\quad \dot{y}(t) = c + u(t),\quad y(t) \le 0,
\]
where $S$ is identified to $\enstq{q \in \R^2}{ 2\, a\, y + b\, x = 0}$ and we assume that $a \ne 0$ and $b \ne 0$. Note that the generalized Legendre-Clebsch condition is related to the sign of $a$. The singular control at $q_0$ is given by $u_s(q_0) = -c-b/2a$ and the boundary control is simply $u_c(q) = -c$. Taking the constraint on the control into account, we see that the condition of admissibility implies
$\abs{c+b/2a}\le 1$ and $\abs{c} \le 1$. Finally, the \emph{hyperbolic case} corresponds to:
\[
a < 0, \quad \abs{c+b/2a} < 1, \quad \abs{c} < 1 \quad \text{and} \quad b\ne 0,
\]
that is the strict generalized Legendre-Clebsch condition is satisfied, the singular and boundary controls are strictly admissible, and the singular set $S$ is not identified with the boundary set $y = 0$.

\begin{figure}[t!]
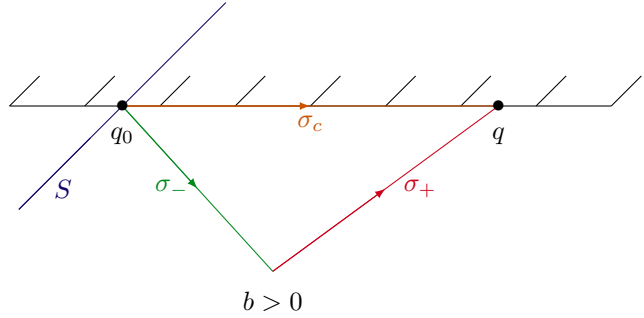

\centering
\boundaryxpos{2.0}
\caption{Small time representation of different trajectories in the neighborhood of a point $q_0$ belonging
to the singular set $S$ and such that $c(q_0)=0$, in the case $b>0$. Two trajectories joining $q_0$ and $q$ are represented:
one boundary arc $\arc_c$ and a bang-bang trajectory $\arc_-\arc_+$. Using the clock form and the Stokes theorem,
we can show that the trajectory $\arc_c$ is time-optimal.}
\label{fig:tps_court}
\end{figure}

For the state unconstrained problem, the singular arc is optimal, each optimal trajectory has at most two switchings and the local optimal synthesis is of the form $\arc_\pm\arc_s\arc_\pm$. For the state constrained problem, the so-called clock form and Stokes theorem (see \cite{bonnard_art:2003,bonnard:2006} for details) can be used to conclude on the optimatility of the boundary arc when $x \ge 0$ and its non-optimality when $x < 0$, for the case $b > 0$. For $b < 0$, the boundary arc is optimal when $x \le 0$ and non-optimal when $x > 0$. See Fig.~\ref{fig:tps_court} for a representation of the bang, singular and boundary arcs in the hyperbolic case with $b > 0$, and for a comparison between boundary and bang-bang arcs.

When $b > 0$, the small time optimal synthesis joining two points $q_1$ and $q_2$ on the boundary $y=0$, with $q_1$ and $q_2$ on each side of $S$ and $q_1$ on the left, is then of the form $\arc_-\arc_s\arc_c$. In this case, each optimal curve in a neighborhood of $q_0$ has at most three switchings and the local optimal synthesis is of the form $\arc_\pm\arc_s\arc_c\arc_\pm$. For $b < 0$, the local optimal synthesis is of the form $\arc_\pm\arc_c\arc_s\arc_\pm$. We summarize this analysis in the following proposition excerpted from \cite{bonnard_art:2003}.

\begin{prpstn}
\label{prop:tps_courts}
Under our assumptions, in the hyperbolic case each small time optimal trajectory has at most three switchings. Moreover, 
\begin{enumerate}
\item For $b<0$, a boundary arc is optimal iff $x\le 0$ and each optimal arc has the form $\arc_\pm\arc_c\arc_s\arc_\pm$.
\item For $b>0$, a boundary arc is optimal iff $x\ge 0$ and each optimal arc has the form $\arc_\pm\arc_s\arc_c\arc_\pm$.
\end{enumerate}
\end{prpstn}

\begin{rmrk}
Prop.~\ref{prop:tps_courts} gives two different patterns $\arc_\pm\arc_s\arc_c\arc_\pm$ and $\arc_\pm\arc_c\arc_s\arc_\pm$,
depending on the local model, for small time optimal trajectories around the boundary.
These two patterns generalize the pattern $\arc_\pm\arc_s\arc_\pm$, presented in Section~\ref{subsec:classification}, in the hyperbolic
state unconstrained case.
\end{rmrk}

\subsection{A first insight on the standard CAS/Mach procedure}
\label{subsec:cas_mach_make_sense}

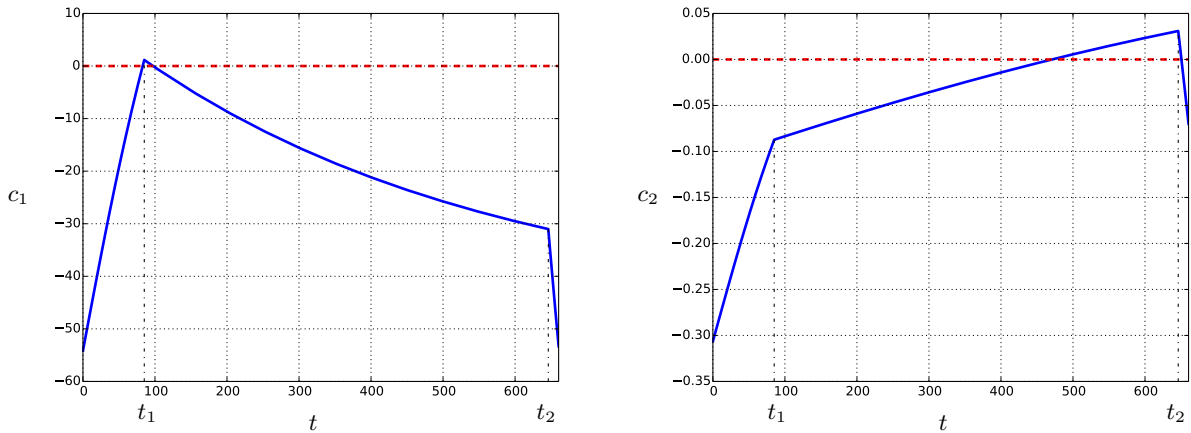
\begin{figure}[t!]
\centering
\def\sizeFig{0.48}
\def\x{203}
\def\y{152}
\begin{tikzgraphics}{\sizeFig\textwidth}{\x}{\y}{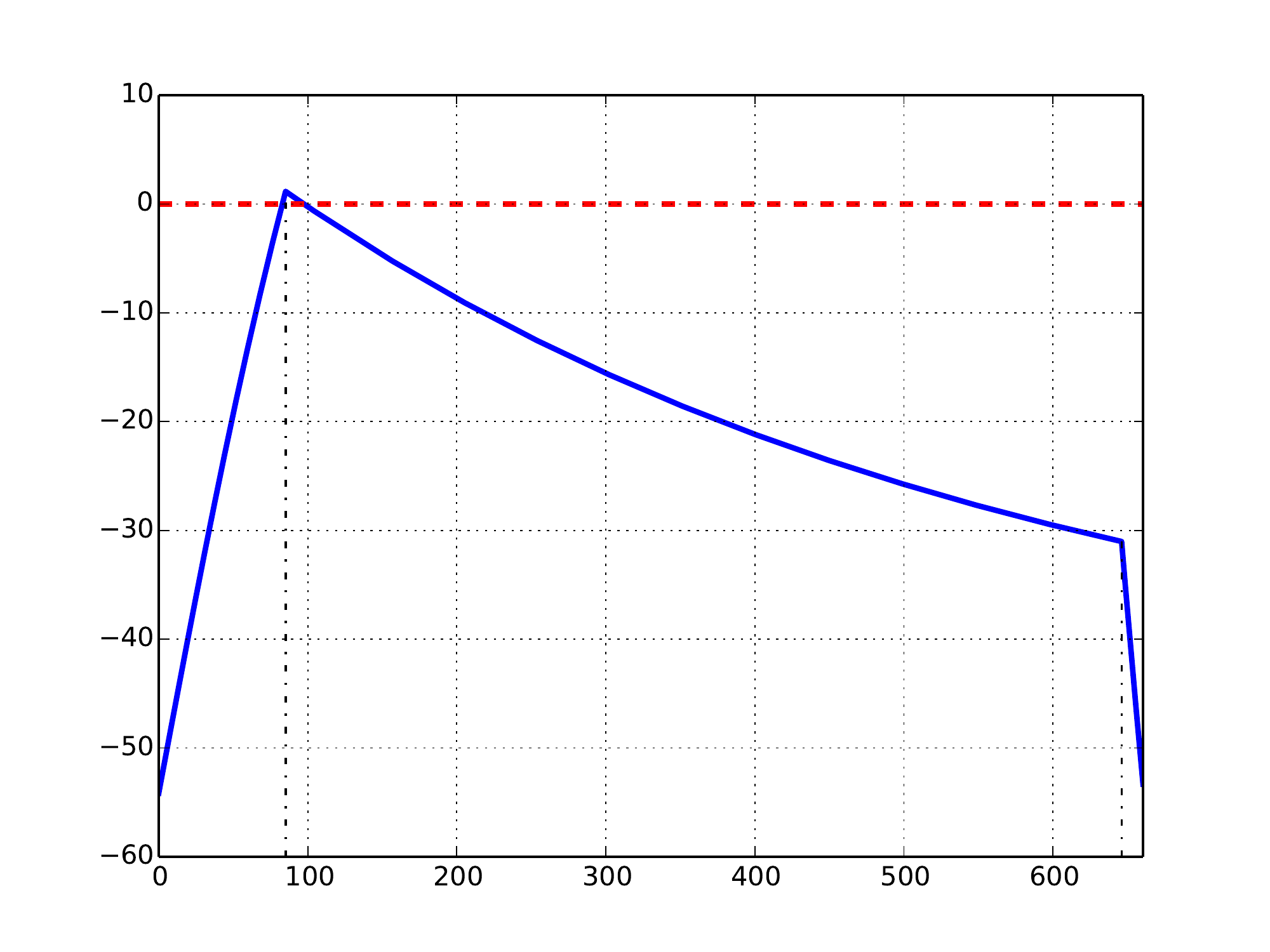}
        \pxcoordinate{0.02*\x}{0.5*\y}{A}; \draw (A) node {\small $c_{1}$};
        \pxcoordinate{0.5*\x}{0.99*\y}{A}; \draw (A) node {\small $t$};
        \pxcoordinate{0.23*\x}{0.97*\y}{A}; \draw (A) node {\small $t_1$};
        \pxcoordinate{0.88*\x}{0.97*\y}{A}; \draw (A) node {\small $t_2$};
\end{tikzgraphics}
\begin{tikzgraphics}{\sizeFig\textwidth}{\x}{\y}{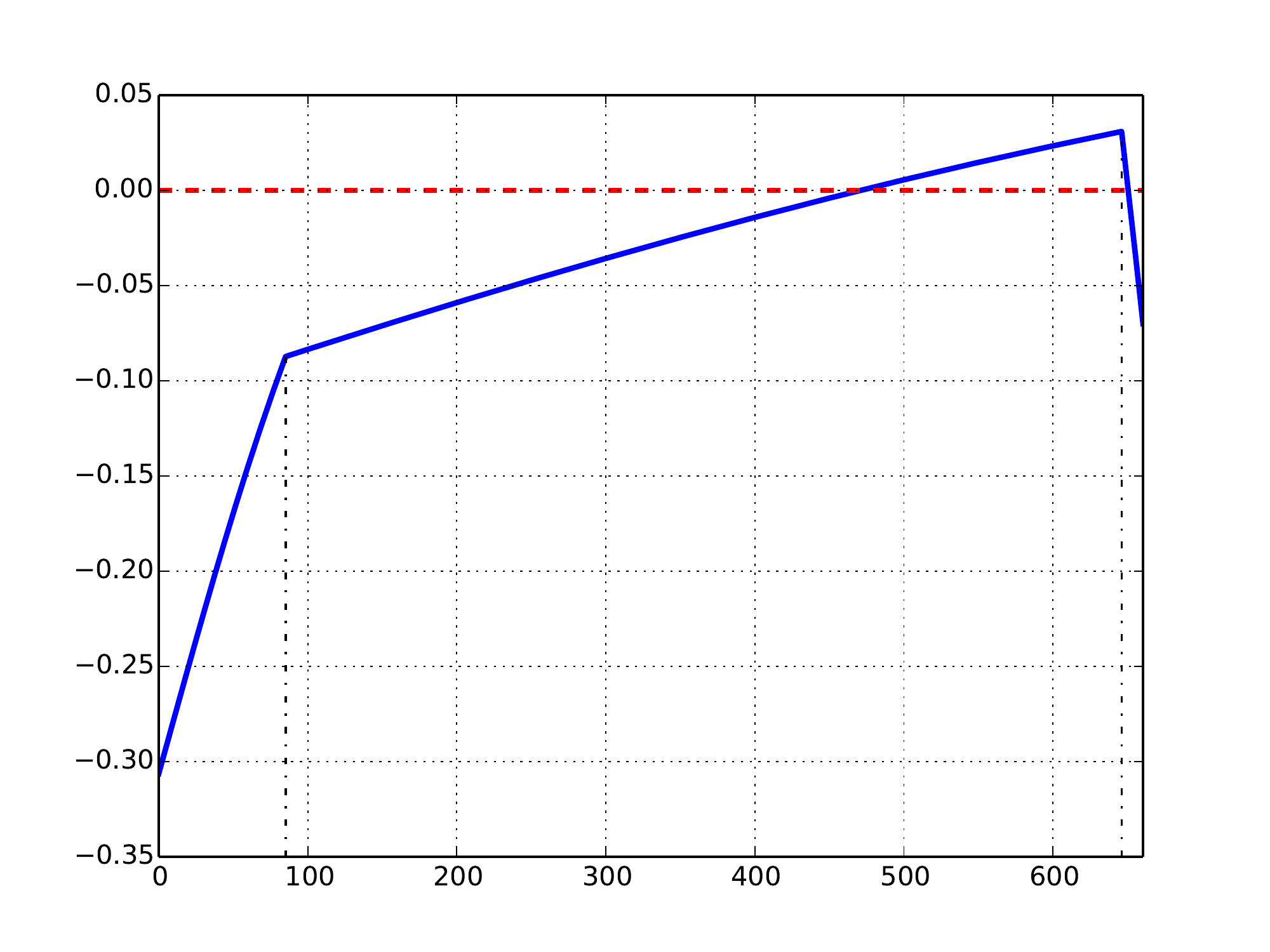}
        \pxcoordinate{0.02*\x}{0.5*\y}{A}; \draw (A) node {\small $c_{2}$};
        \pxcoordinate{0.5*\x}{0.99*\y}{A}; \draw (A) node {\small $t$};
        \pxcoordinate{0.23*\x}{0.97*\y}{A}; \draw (A) node {\small $t_1$};
        \pxcoordinate{0.88*\x}{0.97*\y}{A}; \draw (A) node {\small $t_2$};
\end{tikzgraphics}
\caption{The plain lines depict the evolution of the state constraints $c_1$ and $c_2$ along the state unconstrained hyperbolic trajectory of the form $\sigma_-\sigma_s\sigma_+$, with $\phi_{\max}=160$ and $\psi_{\max}=0.7$. The time $t_1$ (resp. $t_2$) represents the switching time between the negative bang (resp. the singular) arc and the singular (resp. the positive bang) arc.}
\label{fig:struct_constr}
\end{figure}

In this section, we explain on an example why the CAS/Mach procedure makes sense. Let us take from \cite[Fig.~1]{eucass:2017} the state unconstrained trajectory of the form $\arc_-\arc_s\arc_+$, for the minimum time-to-climb problem, with data given in Sections~\ref{sec:intro} and \ref{sec:definitions_ocp}. Note that in \cite{eucass:2017}, $u_\mathrm{min} = -0.262$ instead of $0$ but this has a small influence and it does not change the structure $\arc_-\arc_s\arc_+$. From this excerpted trajectory, denoted $x(\cdot)$, we compute and display on Fig.~\ref{fig:struct_constr}, the state constraints $c_1(x(t)) = \phi(x(t)) - \phimax$ and $c_2(x(t)) = \psi(x(t)) - \psimax$, for $\phimax = 160$ and $\psimax = 0.7$. These realistic values of $\phimax$ and $\psimax$ are chosen to emphasize the following comment: the boundary $c_1 = 0$ is reached during the first bang arc, before the singular arc, while the boundary $c_2 = 0$ is reached during the singular arc,
after reaching $c_1 = 0$ and before the last bang arc.

From Fig.~\ref{fig:struct_constr}, one can also notice that $\phi(x(\cdot))$ and $\psi(x(\cdot))$ reach their maxima respectively at the switching times $t_1$ and $t_2$. Let us consider now other values for $\phimax$ and $\psimax$. Assume we fix
$\phimax = \phi(x(t_1)) - \veps$, $\veps > 0$ small, and $\psimax > \psi(x(t_2))$. Considering these values, the state unconstrained trajectory is not admissible anymore for the state constrained case since in a small neighborhood of $t_1$, we have $c_1(x(t)) > 0$. The trajectory is thus slightly modified and the updated optimal sequence becomes $\arc_-\arc_{c_1}\arc_s\arc_+$ corresponding to the first case of
Prop.~\ref{prop:tps_courts}. On the contrary, if $\phimax > \phi(x(t_1))$ and $\psimax = \psi(x(t_2)) - \veps$, $\veps > 0$ small, then the updated optimal trajectory is of the form $\arc_-\arc_s\arc_{c_2}\arc_+$ corresponding to the second case of Prop.~\ref{prop:tps_courts}. Finally, if we set $\phimax = \phi(x(t_1)) - \veps_1$ and $\psimax = \psi(x(t_2)) - \veps_2$, both $\veps_1 > 0$ and $\veps_2 > 0$ sufficiently small, then the state constrained trajectory is of the form $\arc_-\arc_{c_1}\arc_s\arc_{c_2}\arc_+$, that is the boundary arc $\arc_{c_1}$ appears before the boundary arc $\arc_{c_2}$ and this sequence is made of two patterns. In other words, the trajectory contains one arc at constant CAS which is before the arc at constant Mach, and this is consistent with the actual CAS/Mach procedure. Let us recall that this procedure splits the climbing phase in two parts, the first one is an arc at constant CAS and the second one is an arc at constant Mach number. So, comparing to the CAS/Mach procedure there is an additional singular arc joining the boundary arcs, but, according to the results presented in Section~\ref{sec:minimum_time_application}, for smaller values of $\phimax$ and $\psimax$, the singular arc vanishes and we get the CAS/Mach procedure.

\begin{rmrk}
In Section~\ref{sec:minimum_time_application}, we analyze the influence of $\phimax$ and $\psimax$ on the structure of the trajectories by deforming the BC-extremal associated to the state unconstrained trajectory, taken here as example. Just note that for all the numerical case studies presented in this article, whatever the values of $\phimax$ and $\psimax$, the boundary arc $\arc_{c_1}$ is always before $\arc_{c_2}$.
\end{rmrk}

\section{The minimum time-to-climb problem: numerical methods and results}
\label{sec:minimum_time_application}

\subsection{Introduction to multiple shooting methods}
\label{sec:multiple_shooting}

We present in this section the \emph{indirect multiple shooting method} \cite{bulirsch:2002} that we use to solve the necessary conditions of optimality given
by the maximum principle presented in Section~\ref{sec:general}. We describe the method on only one example and we refer to \cite{theseDamien} for more details. Let us consider the same example as the one presented in the previous section and let us fix
$\phimax > \phi(x(t_1))$ and $\psimax = 0.7$. In this case, the resulting trajectory is of the form  $\arc_-\arc_s\arc_{c_2}\arc_+$.

The unknowns of the shooting method are the junction times $t_1$, $t_2$ and $t_3$,
the final time $t_f$, the jumps $\nu_{t_2}$ and $\nu_{t_3}$ at the junction times $t_2$ and $t_3$ and the initial costate denoted $p_0$. With $p_0$, $t_1$, $t_2$, $t_3$, $\nu_{t_2}$ and $\nu_{t_3}$, since we know the structure of the trajectory,
we can retrieve the states-costates at the junction times $t_1$, $t_2$ and $t_3$ simply by integration, starting from $x(0) = x_0$ and applying the controls $u_\mathrm{min}$, $u_s$, $u_{c_2}$ and $u_\mathrm{max}$, respectively on $\intervalleff{0}{t_1}$, $\intervalleff{t_1}{t_2}$, $\intervalleff{t_2}{t_3}$ and $\intervalleff{t_3}{t_f}$. However, to improve numerical stability, we add the states-costates $z_1$, $z_2$ and $z_3$ at the junction times to the unknowns of the shooting method. The first thing to notice now is that the jumps are zero. Indeed, $\Phi(t_2^-) = 0$ and so $\nu_{t_2} = 0$ (from Prop.~\ref{prop:order1}) since the switching function $\Phi(\cdot) = 0$ along the singular arc $\arc_s$. Now, combining Prop.~\ref{prop:order1} and item 5 of the maximum principle, we have:
\[
\nu_{t_3} = -\frac{\Phi(t_3^+)}{(F_{1} \cdot c_2) (x(t_3))} \le 0.
\]
But, $\Phi(t_3^+) \ge 0$ since the last bang arc is positive. 
{So, assuming $(F_{1} \cdot c_2) (x(t_3)) < 0$, then $\nu_{t_3} = 0$. Note that we can check \textit{a posteriori}
that this condition is verified.}
Finally, grouping all together, we write 
\[
y \coloneqq (p_0, t_1, t_2, t_3, t_f, z_1, z_2, z_3) \in \R^{25}
\]
the unknown of the shooting method.

\begin{rmrk}
For the jump $\nu_{t_3}$ to be not zero, it is necessary that $u_{c_2}(x(t_3)) = u_\mathrm{max}$, which may happen in only few particular cases that we do not encounter in the numerical experiments all through this article.
\end{rmrk}

Now, we need to describe the shooting equations and define the shooting function. First of all, we define the following Hamiltonians:
\begin{equation*}
\begin{aligned}
    H_+(z)      &\coloneqq H_0(z) + u_\mathrm{max}\, H_1(z), \\
    H_-(z)      &\coloneqq H_0(z) + u_\mathrm{min}\, H_1(z), \\
    H_s(z)      &\coloneqq H_0(z) + u_s(x)\, H_1(z), \\
    H_{c_2}(z)  &\coloneqq H_0(z) + u_{c_2}(x)\, H_1(z) + \eta_{c_2}(z)\, c_2(x), \\
\end{aligned}
\end{equation*}
where $z = (x,p)$ and where $u_s$, $u_{c_2}$ and $\eta_{c_2}$ are given in Section~\ref{sec:general}. Note that we can replace $u_s(x) = -D_{001}(x)/D_{101}(x)$ by $u_s(z)=-H_{001}(z)/H_{101}(z)$.
Then, we define the \emph{exponential mapping} $e^{t \varphi}(q_0)$ as the solution at time $t$ of the Cauchy problem $\dot{q}(s) = \varphi(q(s))$, $q(0) = q_0$, where $\varphi$ is any dynamical system and where the state here is denoted $q$. We define also the canonical projections $\pi_x(z) = x$ and $\pi_{p_m}(z) = p_m$, recalling that $z=(x,p)$, $x = (h, v, m)$ and writing $p \coloneqq (p_h, p_v, p_m)$. 
The \emph{multiple shooting function}, denoted $S_{c_2}$, associated to the structure $\arc_-\arc_s\arc_{c_2}\arc_+$ is then defined by:
\[
S_{c_2}(y) \coloneqq 
\begin{pmatrix}
    H_1( e^{t_1 \vv{H_-}}(x_0,p_0) )             \\[0.2em]
    H_{01}( e^{t_1 \vv{H_-}}(x_0,p_0) )          \\[0.2em]
    c_2(\pi_x( e^{(t_2-t_1) \vv{H_s}}(z_1) ))    \\[0.2em]
    \Psi( e^{(t_f-t_3) \vv{H_+}}(z_3) )          \\[0.2em]
    z_1 - e^{t_1 \vv{H_-}}(x_0,p_0)              \\[0.2em]
    z_2 - e^{(t_2-t_1) \vv{H_s}}(z_1)            \\[0.2em]
z_3 - e^{(t_3-t_2) \vv{H_{c_2}}}(z_2)            \\[0.2em]
\end{pmatrix} \in \R^{25},
\]
where 
$\Psi(z) \coloneqq \big(b(\pi_x(z)), \pi_{p_m}(z), H_+(z) + p^0\big)$, $p^0 = -1$,
and the \emph{multiple shooting method} consists in finding a zero of the multiple shooting function $S_{c_2}$, \ie in solving $S_{c_2}(y) = 0$.

\begin{rmrk}
We can replace the first three equations by $H_1(z_1) = H_{01}(z_1) = c_2(\pi_x(z_2)) = 0$.
\end{rmrk}

Let $\ysol \coloneqq (\psol_0, \tsol_1, \tsol_2, \tsol_3, \tsol_f, \zsol_1, \zsol_2, \zsol_3)$ be a zero of the shooting function $S_{c_2}$. Then, to $\ysol$ is associated a unique BC-extremal $(\xsol(\cdot), \psol(\cdot), \usol(\cdot), \etasol(\cdot))$ if
$0 \le \tsol_1 \le \tsol_2 \le \tsol_3 \le \tsol_f$ and if $H_1(\zsol(t)) < 0$ for almost every $t \in \intervalleff{0}{\tsol_1}$ and $H_1(\zsol(t)) > 0$ for almost every $t \in \intervalleff{\tsol_3}{\tsol_f}$.
These two conditions are not included in the shooting equations but they are easily checked \textit{a posteriori} and should be satisfied if the structure is the right one.
Let us explain now why these shooting equations define a BC-extremal of the form $\arc_-\arc_s\arc_{c_2}\arc_+$. The first two equations in $S_{c_2}(y)=0$ impose the extremal to be singular on $\intervalleff{\tsol_1}{\tsol_2}$ according to Prop.~\ref{prop:singular_extremal}, the third equation impose a boundary arc on $\intervalleff{\tsol_2}{\tsol_3}$ according to Prop.~\ref{prop:boundary_extremal}, noticing that $H_1(e^{(\tsol_2-\tsol_1) \vv{H_s}}(\zsol_1))=0$, while the fourth equation $\Psi = 0$ contains the limit and transversality conditions. The last equations are the so-called \emph{matching conditions}.

\begin{figure}[t!]
\centering
\def\sizeFig{0.48}
\def\x{203}
\def\y{152}
\begin{tikzgraphics}{\sizeFig\textwidth}{\x}{\y}{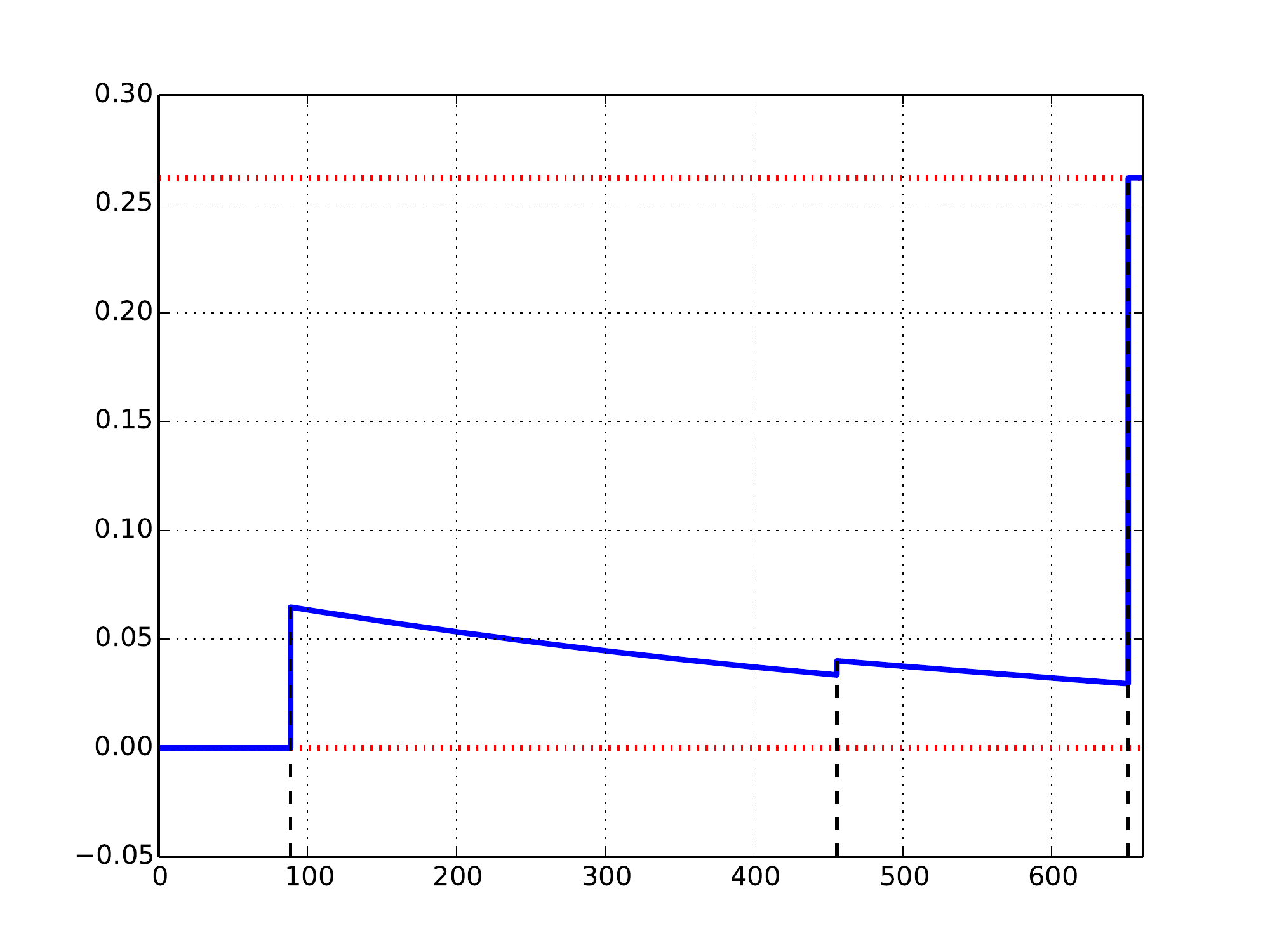}
        \pxcoordinate{0.02*\x}{0.5*\y}{A}; \draw (A) node {\small $u$};
        \pxcoordinate{0.5*\x}{0.99*\y}{A}; \draw (A) node {\small $t$};
        \pxcoordinate{0.23*\x}{0.97*\y}{A}; \draw (A) node {\small $\bar{t}_1$};
        \pxcoordinate{0.67*\x}{0.97*\y}{A}; \draw (A) node {\small $\bar{t}_2$};
        \pxcoordinate{0.89*\x}{0.97*\y}{A}; \draw (A) node {\small $\bar{t}_3$};
\end{tikzgraphics}
\begin{tikzgraphics}{\sizeFig\textwidth}{\x}{\y}{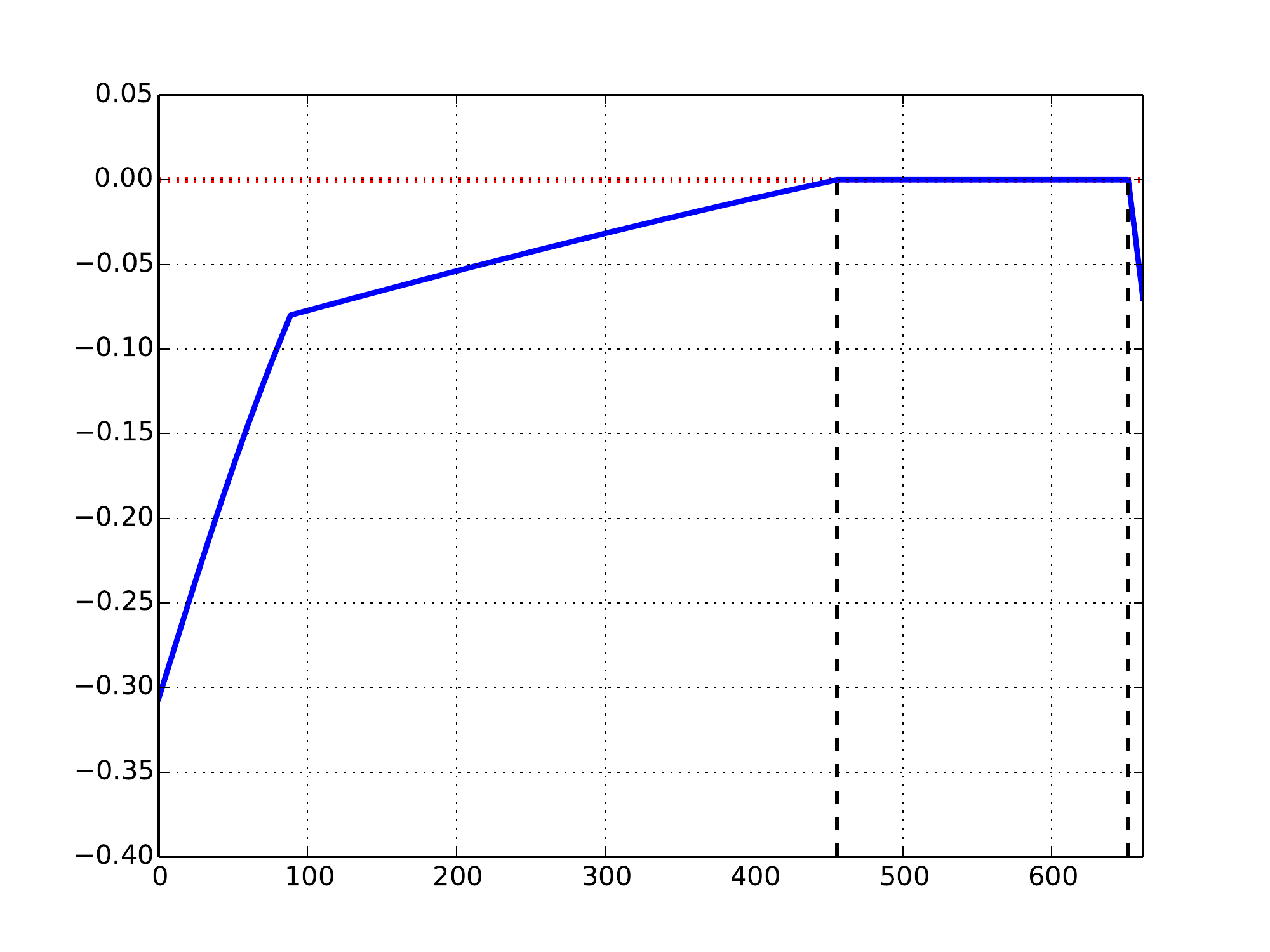}
        \pxcoordinate{0.02*\x}{0.5*\y}{A}; \draw (A) node {\small $c_{2}$};
        \pxcoordinate{0.5*\x}{0.99*\y}{A}; \draw (A) node {\small $t$};
        \pxcoordinate{0.67*\x}{0.97*\y}{A}; \draw (A) node {\small $\bar{t}_2$};
        \pxcoordinate{0.89*\x}{0.97*\y}{A}; \draw (A) node {\small $\bar{t}_3$};
\end{tikzgraphics}
\caption{Evolution of the control and the constraint $c_2$ along the trajectory associated to the zero $\ysol$ of the shooting function $S_{c_2}$.}
\label{fig:sol_contraint}
\end{figure}
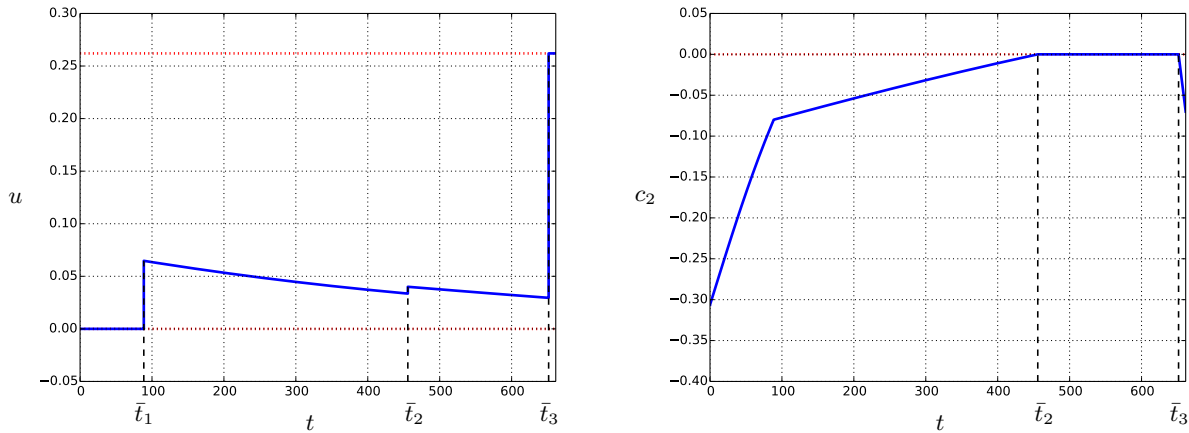

We use the \hampath\ software \cite{hampath:2012, cots:2017}, to obtain a (numerical) zero, denoted $\ysol$, of the shooting function with high accuracy: $\norme{S_{c_2}(\ysol)} \approx 1.511 \times 10^{-11}$. The switching times are $\tsol_1 \approx 88.61$s, $\tsol_2 \approx 455.7$s and $\tsol_3 \approx 651.46$s.
The final time is $\tsol_f \approx 661.37$s and the initial adjoint vector is
$\psol_0 \approx (3.180 \times 10^{-2}, 6.594\times 10^{-1}, -2.300\times 10^{-1})$.
Note that the \hampath\ software is based on a Newton-like algorithm which is initialized thanks to the state unconstrained trajectory. Fig.~\ref{fig:sol_contraint} depicts the control and the constraint $c_2$ along the
resulting state constrained trajectory.

\subsection{Introduction to differential path following methods with monitoring}

The shooting method is used to solve a single optimal control problem. To solve a one-parameter family of optimal control problems, \textit{e.g.} for different values of $\psimax$, we use \emph{differential path following} techniques \cite{allgower:2003,gergaud:2008} with arc length parameterization (also called homotopy method).
Let $h \colon \R^N\times\R\rightarrow\R^N$, $h(y,\lambda)$, denote an \emph{homotopic function} with $\lambda$ the \emph{homotopic parameter}. For example, one can consider the homotopic function defined by $S_{c_2}$, with $y \coloneqq (p_0, t_1, t_2, t_3, t_f, z_1, z_2, z_3)$, $N = 25$ and $\lambda \coloneqq \psimax$. That is to say, in this example $h(y,\lambda) = S_{c_2}(y)$ with $\lambda = \psimax$ which appears in the constraint function $c_2$. We are interested in solving $h = 0$ for $\lambda \in \intervalleff{\lambda_0}{\lambda_1}$, $(\lambda_0, \lambda_1) \in \R^2$ given.

The classical difficulties about homotopic methods consist in assuring that a curve in $h^{-1}(0)$ exists, is sufficiently smooth and will intersect the fixed target $\lambda_1$ in finite length starting from $\lambda=\lambda_0$. As a first result we have the following. Suppose $h$ is continuously differentiable and that we know $y_0$ such that $h(y_0,\lambda_0)=0$ and assume 
\[
\rank \left(\frac{\partial h}{\partial y}(y_0,\lambda_0)\right) = N.
\]
Suppose also that $0$ is a regular value of $h$. Then, a continuously differentiable curve crossing $(y_0, \lambda_0)$ and transverse to $\lambda = \lambda_0$ exists and is either diffeomorphic to a circle or the real line. The curves in $h^{-1}(0)$ are disjoints, and we call each branch of $h^{-1}(0)$ a path of zeros. To compute numerically these paths of zeros, we use the homotopic method from the \hampath\ software, which is based on a Predictor-Corrector algorithm with a high order
and variable step size Runge-Kutta scheme for the prediction and with a classical simplified Newton method for the correction.


We present in this section how we combine homotopy and \emph{monitoring} (that we define hereinafter) to obtain, for a fixed climbing scenario, a cartography of possible time-optimal structures with respect to the values of the bounds of the state constraints $c_1$ and $c_2$, \textit{i.e.} on the values of $\phimax$ and $\psimax$.
We do not give all the details and we refer to \cite{cots:2017} for a more detailed description of the methodology to obtain such a cartography. We start from the state unconstrained trajectory, denoted $x(\cdot)$, presented in Section~\ref{subsec:cas_mach_make_sense} and of the form $\sigma_-\sigma_s\sigma_+$. The maximal values of $\phi$ and $\psi$ along the trajectory are reached respectively at the switching times $t_1$ and $t_2$. One can notice that $\phi(x(t_1)) < VMO$ and $\psi(x(t_2)) < MMO$, so for this trajectory, the state constraints $c_1$ and $c_2$ are not violated if we fix $\phimax = VMO$ and $\psimax = MMO$. The idea is then to deform by homotopy, making $\phimax$ or $\psimax$ varying, the BC-extremal associated to this trajectory and then to detect by monitoring when a change in the structure occurs. The monitoring consists in checking some conditions after each step of prediction and correction. Technically, we only have to code the conditions and then the \hampath\ software automatically check them after each prediction-correction step and stop the homotopy process if at least
one condition is violated. Here are the three different monitorings we use:
\begin{itemize}
\item[\textbf{M1}] check if $c_1$ is violated at the entry point of $\sigma_s$. If yes, then we add a $\arc_{c_1}$ arc;
\item[\textbf{M2}] check if $c_2$ is violated at the exit point of $\sigma_s$. If yes, then we add a $\arc_{c_2}$ arc;
\item[\textbf{M3}] check if $\sigma_s$ has a positive length. If not, then the arc is removed.
\end{itemize}
When a change is detected, we update the homotopic function accordingly to the
new structure and repeat the process. We thus have to limit the range of values for $\phimax$ and $\psimax$. We choose to build the cartography for $(\phimax, \psimax) \in \intervalleff{\phi_0}{VMO} \times \intervalleff{\psi_0}{MMO} $, where $VMO \approx 180$, $MMO \approx 0.82$, $\phi_0 \coloneqq \max(\phi(x(0)), \phi(x(t_f))) \approx 107$ and $\psi_0 \coloneqq \max(\psi(x(0)), \psi(x(t_f))) \approx 0.63$, with $t_f$ the final time. With these chosen minimal values $\psi_0$ and $\phi_0$, whatever $\phimax$ and $\psimax$, the trajectory will always start with an arc $\arc_-$ and end with an arc $\arc_+$ because the state constraints depend only on $h$ and $v$ which are fixed at the initial and final times to some given values which do not depend on $\phimax$ and $\psimax$. This explains why we say that we study the deformation of the BC-extremal for a fixed climbing scenario.

\subsection{Classification of BC-extremals}

Let simply recall that we start from the BC-extremal associated to the state unconstrained trajectory of the form $\arc_-\arc_s\arc_+$, denoted $x(\cdot)$, with $t_1$ and $t_2$ the switching times and $t_f$ the final time.

\medskip
\paragraph{\textbf{Step 1.}}
Let $\phi_{c_1} \coloneqq \phi(x(t_1))$, $\psi_{c_2} \coloneqq \psi(x(t_2))$ and let us label $\delta_{c_1}$, $\delta_{c_2}$ and $\delta_{s}$, a touch point respectively with the state constraint $c_1=0$, $c_2=0$ and the singular set $\Sigma_s$. Then, we have the following straightforward first result:
\begin{center}
\begin{tabular}{ccc}
\medhrule
$\phimax$                           & $\psimax$                         & structure                 \\[0.5em]
\bighrule
    $\phi_{c_1}$
&   $\psi_{c_2}$
&   $\arc_-\delta_{c_1}\arc_s\delta_{c_2}\arc_+$                        \\[0.5em]
    $\phi_{c_1}$
&   $\intervalleof{\psi_{c_2}}{MMO}$
&   $\arc_-\delta_{c_1}\arc_s\arc_+$                                     \\[0.5em]
    $\intervalleof{\phi_{c_1}}{VMO}$
&   $\psi_{c_2}$
&   $\arc_-\arc_s\delta_{c_2}\arc_+$                                     \\[0.5em]
    $\intervalleof{\phi_{c_1}}{VMO}$
&   $\intervalleof{\psi_{c_2}}{MMO}$
&   $\arc_-\arc_s\arc_+$                                                  \\[0.5em]
\medhrule
\end{tabular}
\end{center}

\medskip
\paragraph{\textbf{Step 2.}}
Let fix first $\psimax = MMO$. For $\phimax = \phi_{c_1} - \veps$, $\veps > 0$ small enough, the BC-extremal has a structure of the form $\arc_-\arc_{c_1}\arc_s\arc_+$. Then, we perform an homotopy on $\lambda = \phimax$ from $\lambda_0 = \phi_{c_1}$ to $\lambda_1 = \phi_0$ which stops because of monitoring \textbf{M3} around $\lambda \approx 129.8 \eqqcolon \phi_s$. Hence, at $\lambda = \phi_s$ the structure is $\arc_-\arc_{c_1} \delta_s\arc_+$. The second homotopy from $\lambda_0 = \phi_s$ to $\lambda_1 = \phi_0$ with a structure $\arc_-\arc_{c_1}\arc_+$ is not stopped by the monitoring, thus, no change in the structure is detected anymore. Let fix now $(\phimax, \psimax) = (\phi_s, MMO)$. The structure is of the form $\arc_-\arc_{c_1} \delta_s\arc_+$. Then, we define $\psi_s \coloneqq \max \psi(\cdot) \approx 0.734$ along the corresponding trajectory, and we have the following:
\begin{center}
\begin{tabular}{ccc}
\medhrule
$\phimax$                           & $\psimax$                         & structure                 \\[0.5em]
\bighrule
    $\intervallefo{\phi_0}{\phi_s}$
&   MMO
&   $\arc_-\arc_{c_1}\arc_+$                              \\[0.5em]
    $\phi_s$
&   MMO
&   $\arc_-\arc_{c_1} \delta_s\arc_+$                     \\[0.5em]
    $\intervalleoo{\phi_s}{\phi_{c_1}}$
&   MMO
&   $\arc_-\arc_{c_1}\arc_s\arc_+$                       \\[0.5em]
\medhrule
    $\phi_s$
&   $\psi_s$
&   $\arc_-\arc_{c_1} \delta_s \delta_{c_2}\arc_+$        \\[0.5em]
    $\phi_s$
&   $\intervalleof{\psi_s}{MMO}$
&   $\arc_-\arc_{c_1} \delta_s\arc_+$                     \\[0.5em]
\medhrule
\end{tabular}
\end{center}

\begin{figure}[t!]
\centering
\def\sizeFig{0.45}
\def\x{203}
\def\y{152}
\begin{tikzgraphics}{\sizeFig\textwidth}{\x}{\y}{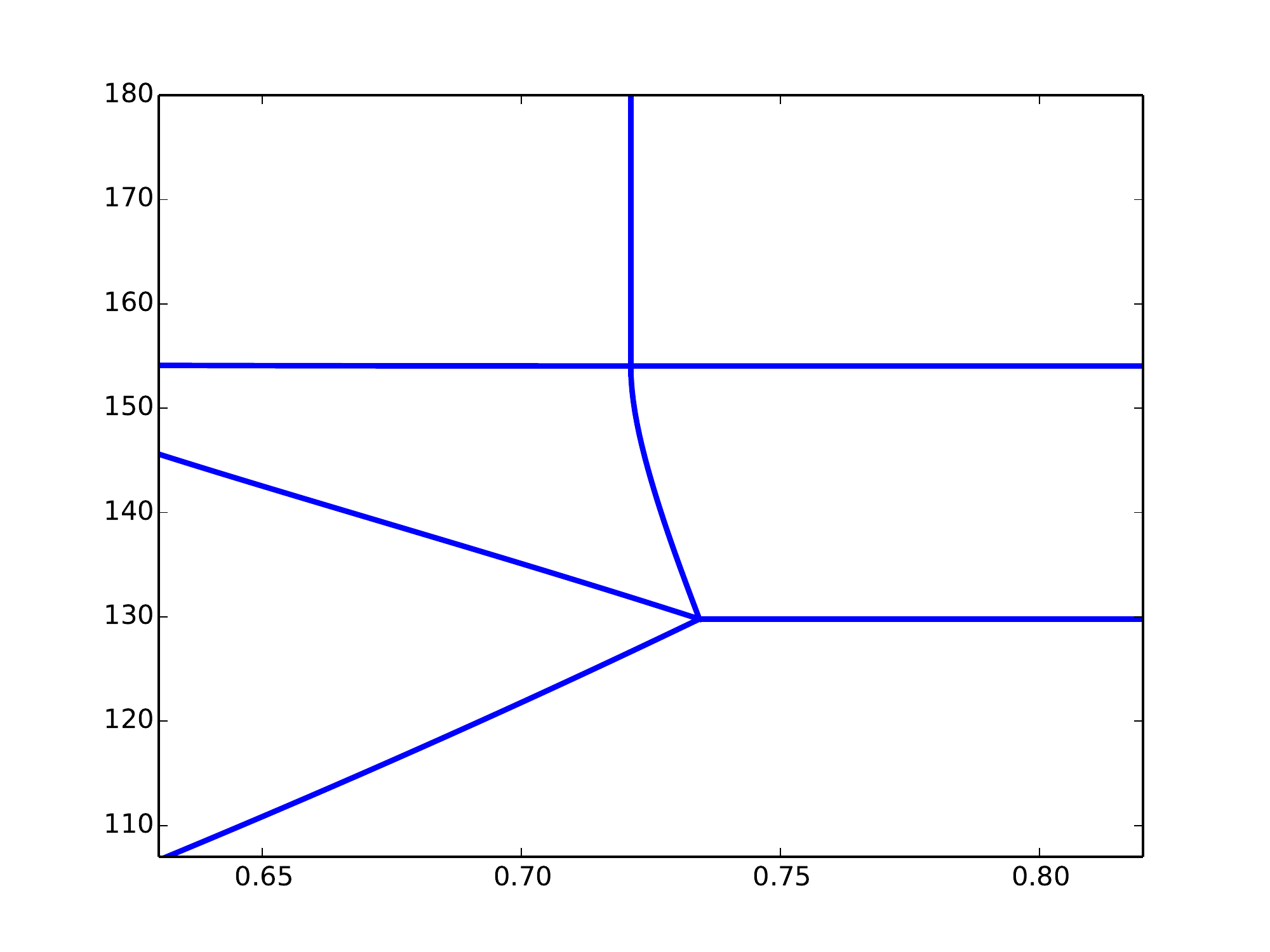}
\pxcoordinate{0.02*\x}{0.5*\y}{A}; \draw (A) node {\small $\phimax$};
\pxcoordinate{0.5*\x}{1.05*\y}{A}; \draw (A) node {\small $\psimax$};
\pxcoordinate{0.497*\x}{0.96*\y}{A}; \draw (A) node {\small $\psi_{c_2}$};
\pxcoordinate{0.56*\x}{0.96*\y}{A}; \draw (A) node {\small $\psi_{s}$};
\pxcoordinate{0.13*\x}{0.96*\y}{A}; \draw (A) node {\small $\psi_{0}$};
\pxcoordinate{0.91*\x}{0.96*\y}{A}; \draw (A) node {\small MMO};
\pxcoordinate{0.95*\x}{0.38*\y}{A}; \draw (A) node {\small $\phi_{c_1}$};
\pxcoordinate{0.97*\x}{0.1*\y}{A}; \draw (A) node {\small VMO};
\pxcoordinate{0.95*\x}{0.65*\y}{A}; \draw (A) node {\small $\phi_{s}$};
\pxcoordinate{0.95*\x}{0.89*\y}{A}; \draw (A) node {\small $\phi_{0}$};
\pxcoordinate{0.7*\x}{0.25*\y}{A}; \draw (A) node {\small $\arc_-\arc_s\arc_+$};
\pxcoordinate{0.7*\x}{0.55*\y}{A}; \draw (A) node {\small $\arc_-\arc_{c_1}\arc_s\arc_+$};
\pxcoordinate{0.7*\x}{0.8*\y}{A}; \draw (A) node {\small $\arc_-\arc_{c_1}\arc_+$};
\pxcoordinate{0.3*\x}{0.25*\y}{A}; \draw (A) node {\small $\arc_-\arc_s\arc_{c_2}\arc_+$};
\pxcoordinate{0.32*\x}{0.47*\y}{A}; \draw (A) node {\small $\arc_-\arc_{c_1}\arc_s\arc_{c_2}\arc_+$};
\pxcoordinate{0.28*\x}{0.65*\y}{A}; \draw (A) node {\small $\arc_-\arc_{c_1}\arc_{c_2}\arc_+$};
\pxcoordinate{0.497*\x}{0.900*\y}{A};
\pxcoordinate{0.497*\x}{0.382*\y}{B};
\draw[densely dotted] (A)--(B) ;
\pxcoordinate{0.55*\x}{0.900*\y}{A};
\pxcoordinate{0.55*\x}{0.648*\y}{B};
\draw[densely dotted] (A)--(B) ;
\end{tikzgraphics}
\def\sizeFig{0.45}
\def\x{203}
\def\y{152}
\begin{tikzgraphics}{\sizeFig\textwidth}{\x}{\y}{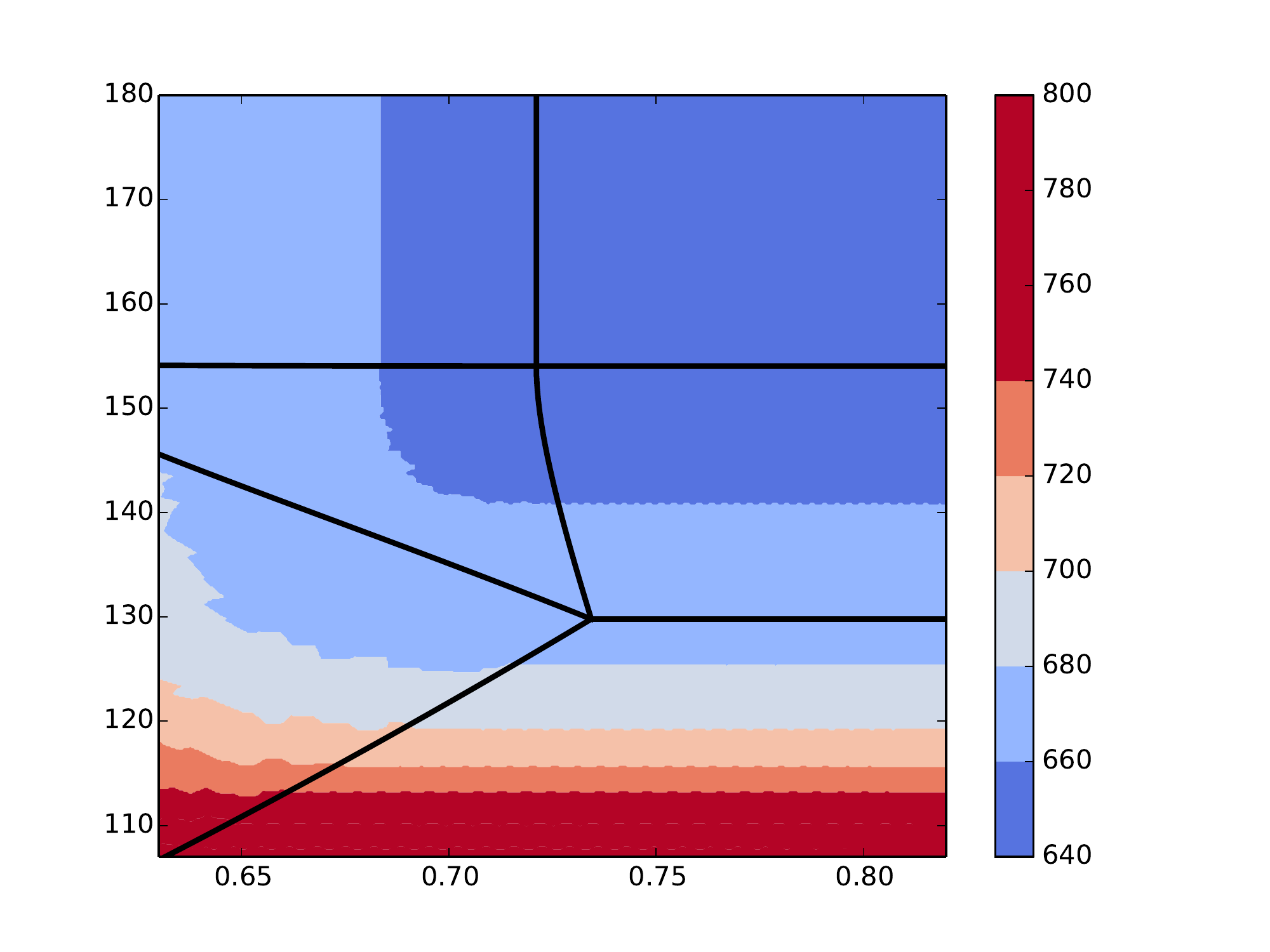}
\pxcoordinate{0.02*\x}{0.5*\y}{A}; \draw (A) node {\small $\phimax$};
\pxcoordinate{0.5*\x}{0.97*\y}{A}; \draw (A) node {\small $\psimax$};
\pxcoordinate{0.5*\x}{1.05*\y}{A}; \draw (A) node {\small \phantom{$\psimax$}};
\end{tikzgraphics}
\caption{(Left) Classification of the BC-extremal structures with respect to $\phimax$ and $\psimax$. To the separating blue lines is associated structures with touch points. These blue lines are computed also by homotopy and monitoring. (Right) The final time with respect to $\phimax$ and $\psimax$.}
\label{fig:classification}
\end{figure}
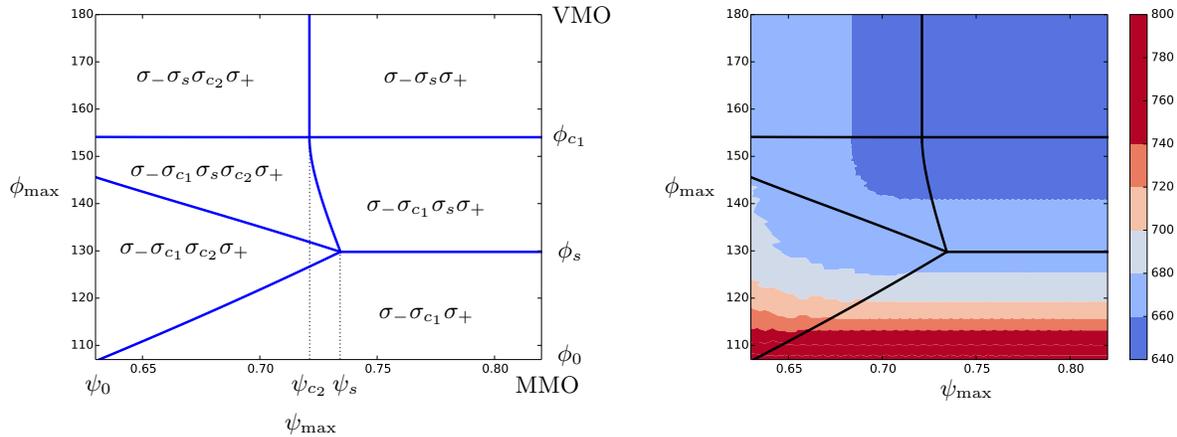

\paragraph{\textbf{Step 3.}}
The same methodology is used each time the homotopy is stopped by a monitoring. The new structure is then guessed thanks to the results from Section \ref{sec:time-to-climb}. The final classification is given on the left sub-graph of Fig.~\ref{fig:classification}.

\subsection{Discussion on the classification}

The left sub-graph of Fig.~\ref{fig:classification} presents the resulting cartography of the BC-extremal structures for $(\phimax, \psimax) \in \intervalleff{\phi_0}{VMO} \times \intervalleff{\psi_0}{MMO}$, for the minimum time-to-climb problem.
The right sub-graph of Fig.~\ref{fig:classification} gives the time-to-climb related to the trajectories from left sub-graph. Of course, the state unconstrained trajectory of the form $\arc_-\arc_s\arc_+$ minimizes the time-to-climb and the minimal value is $t_{f,\min}\coloneqq 658$ s. Note that as a remarkable fact, we retrieve the CAS/Mach procedure for some values of $\phimax$ and $\psimax$.
Up to this point, we can say that the CAS/Mach procedure is not the optimal structure for the minimum time-to-climb problem with the Maximum Operating Speed (VMO) and the Maximum Operating Mach (MMO) speeds as constraints. However, when decreasing these constraints, the CAS/Mach procedure may be optimal since the classification from Fig.~\ref{fig:classification} gives in this case BC-extremals of the form $\arc_-\arc_{c_1}\arc_{c_2}\arc_+$.

\begin{rmrk}
In the previous paragraph, we say ''may be optimal'' since we only check necessary conditions of optimality given by the maximum principle, that is we only compute BC-extremals. The question about local or global optimality is difficult in this context and go beyond the scope of this article.
\end{rmrk}

\begin{rmrk}
According to Section~\ref{sec:minimal_fuel}, for the solutions from Fig.~\ref{fig:classification}, the fuel consumption is minimal for $(\phimax, \psimax) \approx (128.9, 0.661)$ where the structure follows the CAS/Mach procedure.
\end{rmrk}

\section{The CAS/Mach versus the Singular Arc procedure}
\label{sec:cost_index_and_cas_mach}


\subsection{The CAS/Mach (CM) procedure}

This procedure usually splits the climbing phase in two parts, the first one is an arc $\arc_{c_1}$ at constant CAS and the second one is an arc $\arc_{c_2}$ at constant Mach number. In our case, we have seen that we need a first negative bang arc $\arc_{-}$ and we need a final positive bang arc $\arc_{+}$ to reach the limit conditions ($b(x(t_f))=0$) which define the cruise phase. We consider the cost index $g_\alpha$ (see Section~\ref{sec:definitions_ocp}) as the objective function and we want to optimize the CAS and Mach numbers, that is $\phimax$ and $\psimax$, followed respectively along the arcs $\arc_{c_1}$ and $\arc_{c_2}$. The corresponding finite dimensional optimization problem may be summarized this way:
\leqnomode
\begin{equation}
\tag{$\mathrm{CM}_\alpha$}
\min_{X \in \R^{15}} f^{\mathrm{CM}}_{\alpha}(X) \quad \text{subject to} \quad g^{\mathrm{CM}}_{\alpha}(X) = 0,
\label{eq:CM_cost_index}
\end{equation}
\reqnomode
where the decision variable is
$X = (t_1, t_2, t_3, t_f, x_1, x_2, x_3, \phimax, \psimax)$,
where the cost function is given by
\begin{equation*}
\begin{aligned}
    f^{\mathrm{CM}}_{\alpha}(X) &\coloneqq \alpha\, t_f + (1-\alpha) \left[ m_0 - \pi_m\left( e^{(t_f-t_3) F_+}(x_3) \right) \right] \\
    &= g_\alpha\left(t_f, \pi_m\left( e^{(t_f-t_3) F_+}(x_3) \right)\right),
\end{aligned}
\end{equation*}
with $\pi_m(x) = m$, recalling that $x = (h,v,m)$, and where $m_0$ is the initial mass, which is fixed. The equality constraints are defined by:
\[
g^{\mathrm{CM}}_{\alpha}(X) \coloneqq
\begin{pmatrix}
    b( e^{(t_f-t_3) F_+}(x_3) )         \\[0.2em]
    x_1 - e^{t_1 F_-}(x_0)              \\[0.2em]
    x_2 - e^{(t_2-t_1) F_{c_1}}(x_1)    \\[0.2em]
    x_3 - e^{(t_3-t_2) F_{c_2}}(x_2)    \\[0.2em]
    \phimax - \phi(x_1)                 \\[0.2em]
    \psimax - \psi(x_2)                 \\[0.2em]
\end{pmatrix} \in \R^{13},
\]
recalling that $x_0 = (h_0, v_0, m_0)$ is fixed (see Section~\ref{sec:definitions_ocp} for the numerical values).

\begin{rmrk}
One can check \textit{a posteriori} that $0 \le t_1 \le t_2 \le t_3 \le t_f$ or include these inequality constraints into the optimization problem formulation.
\end{rmrk}

\subsection{The Singular Arc (SA) procedure}

The Singular Arc (SA) procedure, which is, according to the next section, better than the CM procedure, is defined by trajectories of the form $\arc_-\arc_s\arc_+$. The cost function is still the cost index and the equality constraints are simply given by the limit and matching conditions. This gives the following NLP problem:
\leqnomode
\begin{equation}
\tag{$\mathrm{SA}_\alpha$}
\min_{X \in \R^{9}} f^{\mathrm{SA}}_{\alpha}(X) \quad \text{subject to} \quad g^{\mathrm{SA}}_{\alpha}(X) = 0,
\label{eq:SA_cost_index}
\end{equation}
\reqnomode
where the decision variable is here
$X = (t_1, t_2, t_f, x_1, x_2)$,
where the cost function is given by
\begin{equation*}
\begin{aligned}
    f^{\mathrm{SA}}_{\alpha}(X) &\coloneqq \alpha\, t_f + (1-\alpha) \left[ m_0 - \pi_m\left( e^{(t_f-t_2) F_+}(x_2) \right) \right] \\
&= g_\alpha\left(t_f, \pi_m\left( e^{(t_f-t_2) F_+}(x_2) \right)\right),
\end{aligned}
\end{equation*}
and where the equality constraints are defined by:
\[
g^{\mathrm{SA}}_{\alpha}(X) \coloneqq
\begin{pmatrix}
    b( e^{(t_f-t_2) F_+}(x_2) )         \\[0.2em]
    x_1 - e^{t_1 F_-}(x_0)              \\[0.2em]
    x_2 - e^{(t_2-t_1) F_{s}}(x_1)      \\[0.2em]
\end{pmatrix} \in \R^{8}.
\]

\begin{rmrk}
One can also check \textit{a posteriori} that $0 \le t_1 \le t_2 \le t_f$  or include these inequality constraints into the optimization problem formulation.
\end{rmrk}

\begin{rmrk}
Let us denote by ($\mathrm{SA}_{t_f}$) the problem \eqref{eq:SA_cost_index} with $\alpha = 1$ and ($\mathrm{CM}_{t_f}$) the problem
\eqref{eq:CM_cost_index} with $\alpha = 1$. According to the results from Section~\ref{sec:minimum_time_application} (see Fig.~\ref{fig:classification} in particular), we have that
($\mathrm{SA}_{t_f}$) is better than ($\mathrm{CM}_{t_f}$).
\end{rmrk}

\subsection{The minimal fuel consumption problem}
\label{sec:minimal_fuel}

We have already compared \eqref{eq:SA_cost_index} and \eqref{eq:CM_cost_index} when $\alpha = 1$. Before comparing them for all possible values of $\alpha$, we present a brief analysis when $\alpha = 0$, that is for the minimum fuel consumption problem.
We denote by ($\mathrm{SA}_{\Delta_m}$) the problem \eqref{eq:SA_cost_index} with $\alpha = 0$ and ($\mathrm{CM}_{\Delta_m}$) the problem
\eqref{eq:CM_cost_index} with $\alpha = 0$.

\begin{figure}[t!]
\centering
\def\sizeFig{0.45}
\def\x{203}
\def\y{152}
\begin{tikzgraphics}{\sizeFig\textwidth}{\x}{\y}{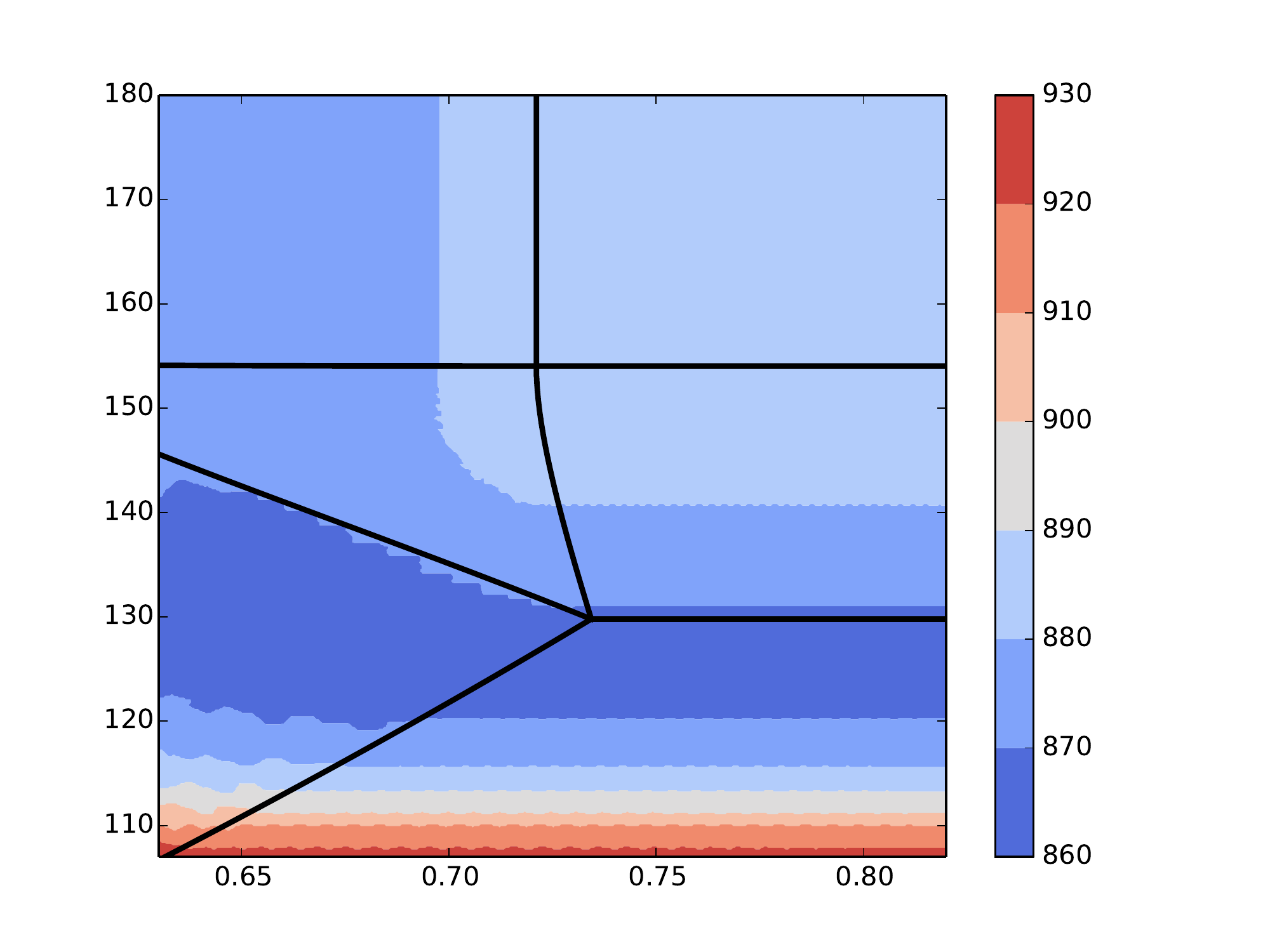}
        \pxcoordinate{0.02*\x}{0.5*\y}{A}; \draw (A) node {\small $\phimax$};
        \pxcoordinate{0.5*\x}{0.97*\y}{A}; \draw (A) node {\small $\psimax$};
\end{tikzgraphics}
\caption{The fuel consumption with respect to $\phimax$ and $\psimax$ for the minimum time-to-climb problem. See the classification on Fig.~\ref{fig:classification}.}
\label{fig:results_homotopy_conso}
\end{figure}

\medskip
\paragraph{\textbf{Fuel consumptions for ($\mathrm{SA}_{t_f}$) versus ($\mathrm{CM}_{\Delta_m}$).}}
For the minimum time-to-climb problem, the classification of the different structures is given on the left sub-graph of Fig.~\ref{fig:classification}. Associated to these BC-extremals, we have the final time $t_f$ presented on the right sub-graph of Fig.~\ref{fig:classification}.
One can notice that the minimal value ($t_{f,\min} = 658$s) is of course for the structure following the SA procedure, and the associated fuel consumption is $\Delta_m \approx 882$kg. For the BC-extremals from Fig.~\ref{fig:classification}, one can also compute the fuel consumption, which is given on Fig.~\ref{fig:results_homotopy_conso}.
Even if it is not clear from the figure, the minimal value of fuel consumption is obtained when $(\phimax, \psimax) \approx (128.9, 0.661)$, that is for a structure following the CM procedure. To obtain these optimal values we solve the problem ($\mathrm{CM}_{\Delta_m}$) with the \textit{fmincon} function from the \matlab\ software and we use the \textit{ode45} function to compute the numerical integrations.
The numerical solution we obtain satisfies the first-order optimality conditions with a tolerance of $10^{-7}$. The switching times are $t_1 \approx 32$s, $t_2 \approx 450$s and $t_3 \approx 676$s while we have a fuel consumption $\Delta_m \approx 863$kg, a final time $t_f \approx 677$s and the CAS and Mach optimal values are $(\phimax, \psimax) \approx (128.9, 0.661)$.
Hence, there is a gap of $19$s on the final time and $19$kg on the fuel consumption between the solution of ($\mathrm{SA}_{t_f}$) and ($\mathrm{CM}_{\Delta_m}$). That is, the decrease of the bounds $\phimax$ and  $\psimax$ of the state constraints
$c_1$ and $c_2$ allows to save up 19 kg of fuel with a 19 seconds longer flight.

\medskip
\paragraph{\textbf{Optimal fuel consumption.}}
In the previous paragraph, we have compared ($\mathrm{SA}_{t_f}$) with ($\mathrm{CM}_{\Delta_m}$). As expected, solving ($\mathrm{SA}_{t_f}$) gives a better final time and solving ($\mathrm{CM}_{\Delta_m}$) gives a better fuel consumption. But we can do better.
Let us come back to the original problem of minimizing the fuel consumption, that is to the optimal control problem ($P_{\Delta_m}$) defined in Section~\ref{sec:definitions_ocp}. To obtain the best possible fuel consumption, we set $\phimax = VMO$ and $\psimax = MMO$.
In this case, thanks to the direct collocation method from the \bocop\ software \cite{bocop:2016}, we determine that the structure $\sigma_-\sigma_s\sigma_+$ is still relevant, that is the solution we obtain follows the SA procedure. This means that the problem ($\mathrm{SA}_{\Delta_m}$) gives a better solution than the problem ($\mathrm{CM}_{\Delta_m}$). 
To get a BC-extremal solution of ($P_{\Delta_m}$) with $\phimax = VMO$ and $\psimax = MMO$, we use the solution from the \bocop\ software to initialize the shooting method. This technique is classical, see \cite{Bonnard2015,theseDamien} for details and examples. 
We then obtain a trajectory for which the switching times are $t_1 \approx 47$s and $t_2 \approx 668$s, the final time is $t_f \approx 675$s and the fuel consumption is $\Delta_m \approx 860$kg. There is a gap of $3$kg on the fuel consumption with the solution from ($\mathrm{CM}_{\Delta_m}$) and a difference about $2$s on the final time. Hence, even if ($\mathrm{SA}_{\Delta_m}$) is better than ($\mathrm{CM}_{\Delta_m}$), the difference is very small.

\begin{rmrk}
The \bocop\ software transforms an infinite dimensional optimal control problem into a finite dimensional optimization problem, applying a full time discretization of the state variables, the control variables, the dynamics and the constraints.
\end{rmrk}

Let us recap in the following table, the fuel consumptions and the final times for the different problems we have encountered so far:
\begin{center}
    \begin{tabular}{c||c|c||c|c}
        \medhrule
            & ($\mathrm{SA}_{t_f}$)  & ($\mathrm{CM}_{t_f}$) & 
            ($\mathrm{SA}_{\Delta_m}$) & ($\mathrm{CM}_{\Delta_m}$)    \\[0.5em]
        \bighrule
            $t_f$       & 658s      & 660s      & 675s      & 677s  \\[0.5em]
            $\Delta_m$  & 882kg     & 884kg     & 860kg     & 863kg \\[0.5em]
        \medhrule
    \end{tabular}
\end{center}
From this table, it is clear that the SA procedure is slightly better than the CM procedure in terms of time-to-climb and fuel consumption.
In the next section we analyze the case when $\alpha \in \intervalleff{0}{1}$. This will complete the comparison between \eqref{eq:CM_cost_index} and
\eqref{eq:SA_cost_index}.

\subsection{The minimal cost index problem}
\label{sec:comparison_cas_mach}

We are interested now in the so-called cost index criterion which is a convex combination of the time-to-climb and the fuel consumption. Instead of comparing directly \eqref{eq:CM_cost_index} and \eqref{eq:SA_cost_index} for $\alpha \in \intervalleff{0}{1}$, we solve \eqref{eq:ocp_cost_index} with $\phimax = VMO$, $\psimax = MMO$, and we show that the trajectories we get follow the SA procedure.
Doing this, we conclude at the end of this section, that not only the SA procedure is better than the CM procedure but it is also the best solution we have. The problem \eqref{eq:CM_cost_index} is solved with the \textit{fmincon} function from the \matlab\ software and we use the \textit{ode45} function to compute the numerical integrations. The optimal control problem \eqref{eq:ocp_cost_index} is solved with the \bocop\ software.

\begin{figure}[t!]
\centering
\def\sizeFig{0.4}
\def\x{167}
\def\y{132}
\begin{tikzgraphics}{\sizeFig\textwidth}{\x}{\y}{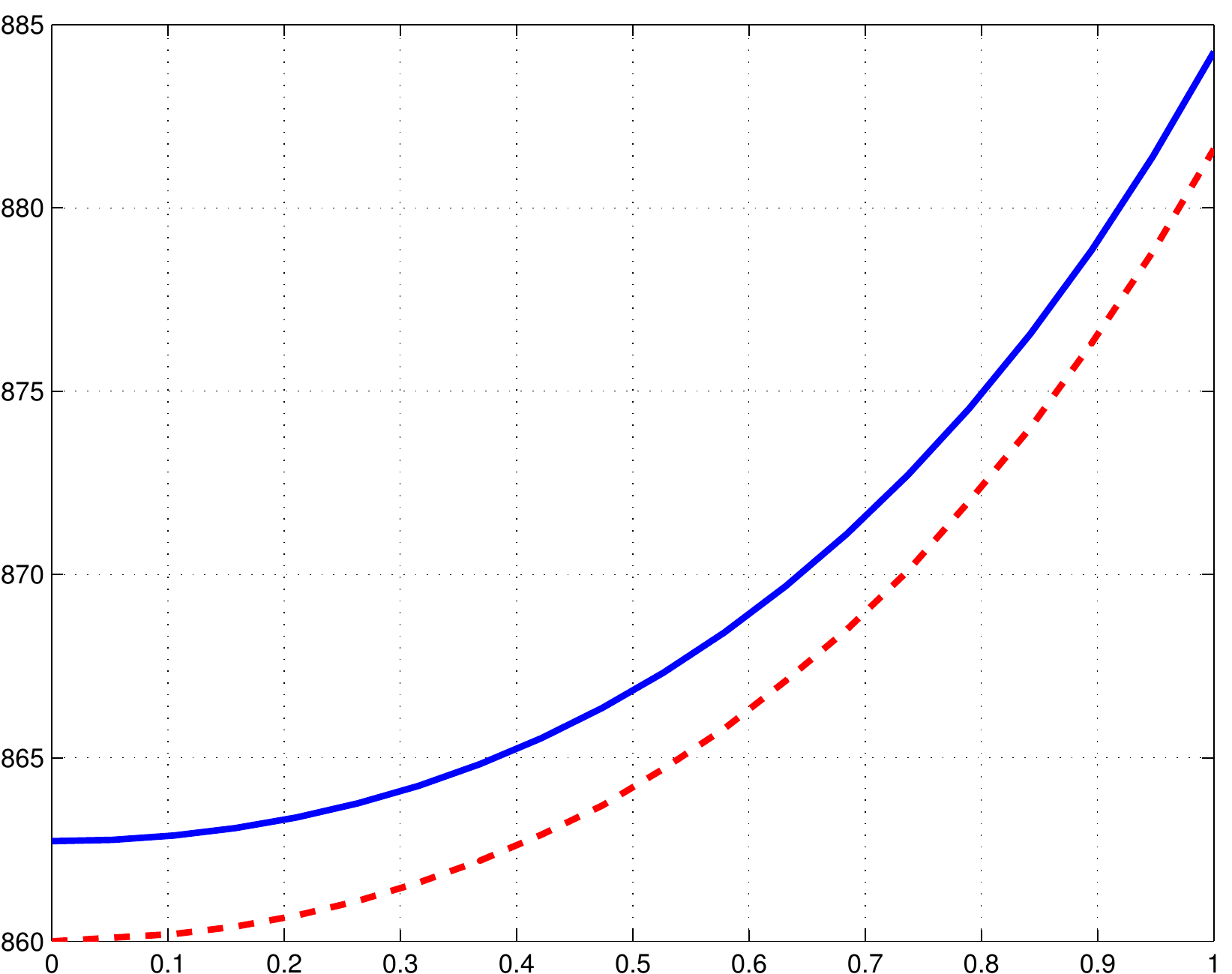}
	\pxcoordinate{-0.03*\x}{0.5*\y}{A}; \draw (A) node {\small $\Delta_m$};
	\pxcoordinate{0.5*\x}{1.03*\y}{A}; \draw (A) node {\small $\alpha$};
\end{tikzgraphics}
\hspace{1cm}
\begin{tikzgraphics}{\sizeFig\textwidth}{\x}{\y}{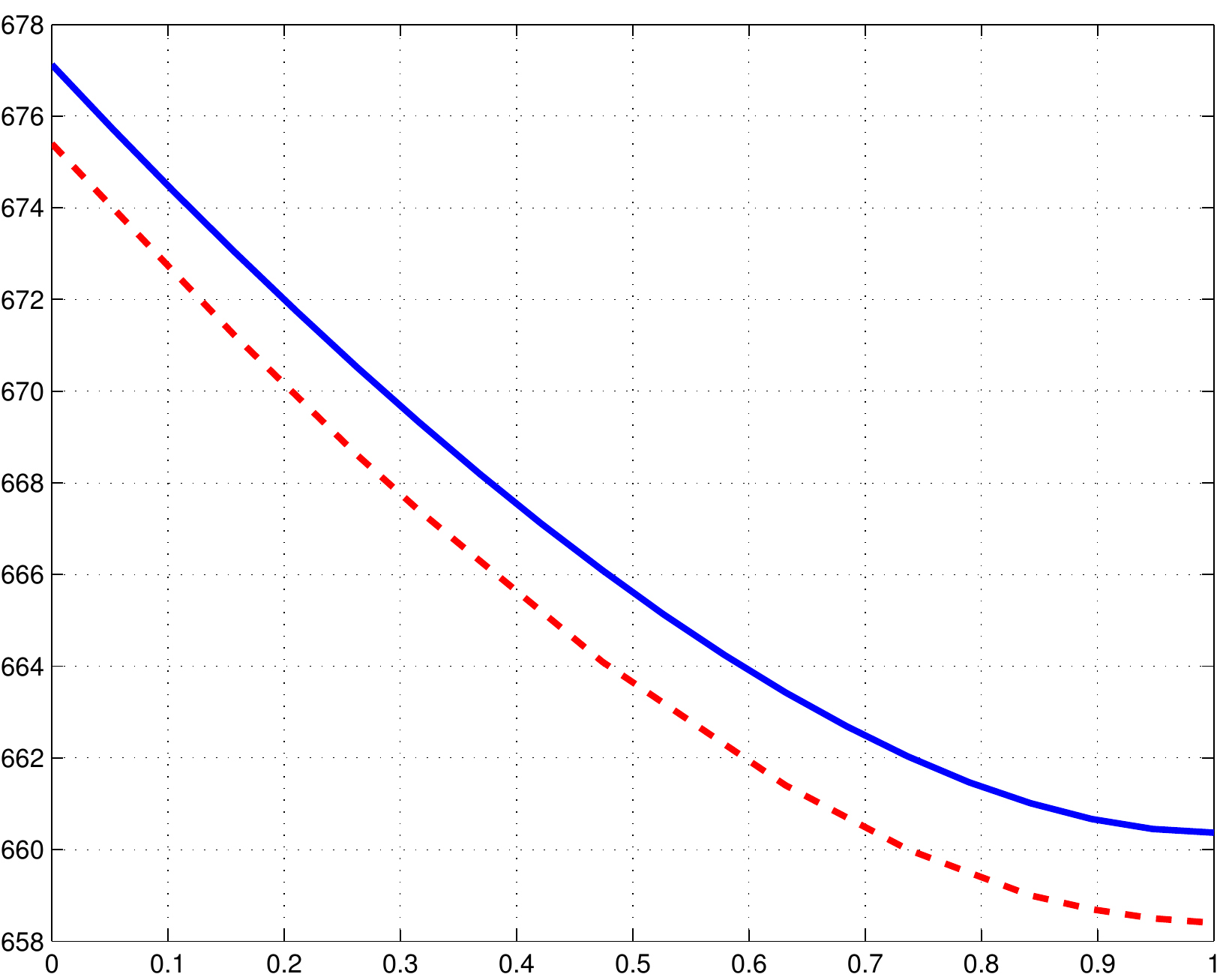}
	\pxcoordinate{-0.04*\x}{0.5*\y}{A}; \draw (A) node {\small $t_f$};
	\pxcoordinate{0.5*\x}{1.03*\y}{A}; \draw (A) node {\small $\alpha$};
\end{tikzgraphics}
\caption{Evolution of the fuel consumption (on the left) and of the time-to-climb (on the right) with respect to $\alpha$. The plain lines are for problem \eqref{eq:CM_cost_index} while the dotted lines are for  problem \eqref{eq:ocp_cost_index} with $\phimax = VMO$, $\psimax = MMO$. The two curves are fairly similar, the relative gap between the plain and the dotted curves is about $0.3\%$.}
\label{fig:comp_nlp_ocp}
\end{figure}

We present in Table~\ref{tab:nlp_output}, for different values of $\alpha$, the optimized CAS/Mach couple (that is $\phimax$ and $\psimax$), the fuel consumption and the time-to-climb, associated to the solutions of problems \eqref{eq:CM_cost_index} and \eqref{eq:ocp_cost_index}.
Since all the trajectories from problem \eqref{eq:ocp_cost_index} with $\phimax = VMO$, $\psimax = MMO$ follow the SA procedure, it is clear from Table~\ref{tab:nlp_output} that the SA procedure is better than the CM procedure.
Indeed, for any $\alpha$, the final time and the fuel consumption are better, so is the cost index. It is interesting to notice that not only the objective function is better (that is the cost index) but also each part of the criterion (that is $t_f$ and $\Delta_m$). 
This fact is visible on Fig.~\ref{fig:comp_nlp_ocp}, where the left sub-graph depicts the evolution of the fuel consumption with respect to $\alpha$ while the right sub-graph presents the time-to-climb.

To conclude, this numerical investigation indicates that the SA procedure of the form $\arc_-\arc_s\arc_+$ is the best policy but the actual CAS/Mach procedure of the form $\arc_-\arc_{c_1}\arc_{c_2}\arc_+$ provides good sub-optimal trajectories regarding the cost index criterion $g_\alpha$.

\begin{table}[t!]
\centering
\begin{tabular}{c||c|c|c|c||c|c}
 &\multicolumn{4}{c||}{\eqref{eq:CM_cost_index}} & \multicolumn{2}{c}{\eqref{eq:ocp_cost_index}}\\
$\alpha$&CAS & Mach & $\Delta_{m}$ & $t_f$ & $\Delta_m$ & $t_f$\\
\hline \hline
0     & 128.9 & 0.6611 & 862.7  & 677.1        & 860.0    & 675.4 \\   
0.056 & 129.4 & 0.6631 & 862.8  & 675.7        & 860.1    & 674.0 \\   
0.105 & 129.9 & 0.6651 & 862.9  & 674.3        & 860.2    & 672.6 \\   
0.158 & 130.4 & 0.6670 & 863.1  & 673.0        & 860.4    & 671.2 \\   
0.210 & 131.0 & 0.6690 & 863.4  & 671.7        & 860.7    & 669.9 \\
0.263 & 131.6 & 0.6709 & 863.8  & 670.5        & 861.1    & 668.6 \\
0.316 & 132.2 & 0.6729 & 864.2  & 669.3        & 861.6    & 667.4 \\
0.368 & 132.8 & 0.6748 & 864.8  & 668.2        & 862.2    & 666.3 \\
0.421 & 133.5 & 0.6767 & 865.5  & 667.1        & 862.9    & 665.2 \\
0.474 & 134.2 & 0.6787 & 866.4  & 666.1        & 863.7    & 664.1 \\
0.526 & 134.9 & 0.6808 & 867.3  & 665.1        & 864.7      & 663.2 \\   
0.579 & 135.6 & 0.6825 & 868.4  & 664.2        & 865.8      & 662.3 \\   
0.631 & 136.4 & 0.6844 & 869.7  & 663.4        & 867.1      & 661.4 \\   
0.684 & 137.2 & 0.6863 & 871.1  & 662.7        & 868.5      & 660.7 \\   
0.737 & 138.1 & 0.6882 & 872.7  & 662.0        & 870.1      & 660.0 \\   
0.790 & 139.0 & 0.6901 & 874.5  & 661.5        & 872.0      & 659.5 \\   
0.842 & 139.9 & 0.6920 & 876.6  & 661.0        & 874.0      & 659.0 \\   
0.895 & 140.9 & 0.6939 & 878.8  & 660.7        & 876.3      & 658.7 \\   
0.947 & 142.0 & 0.6959 & 881.4  & 660.4        & 878.8      & 658.5 \\   
1.000 & 143.2 & 0.6978 & 884.3  & 660.4        & 881.6      & 658.4    
\end{tabular}
\caption{Fuel consumption and time-to-climb for different values of $\alpha$ associated to the solutions of problems \eqref{eq:CM_cost_index}
and \eqref{eq:ocp_cost_index}. For problem \eqref{eq:ocp_cost_index}, we set $\phimax = VMO$ and $\psimax=MMO$.}
\label{tab:nlp_output}
\end{table}

\section{Conclusion}
\label{sec:conclusion}

This article is about aircraft trajectory optimization during the climbing phase considering the cost index criterion. The cost index criterion is a convex combination of the time-to-climb and the fuel consumption, parameterized by $\alpha$. This parameter $\alpha$ is tuned by each airline.
The problem is first modeled as an optimal control problem in Mayer form with a single-input affine control system and with pure state constraints. The two state constraints give CAS and Mach speeds limitations while the control variable is the {flight path angle} of the aircraft. In Section~\ref{sec:minimum_time_application}, for the minimum time-to-climb problem ($\alpha = 0$), we classify the BC-extremals structures with respect to the bounds of the state constraints, that is $\phimax$ and $\psimax$, using the small time analysis developed in Section~\ref{sec:time-to-climb}, combined with indirect multiple shooting and homotopy methods with monitoring. %
This classification emphasizes the role of the SA and CM procedures that we define and compare in Section~\ref{sec:cost_index_and_cas_mach}. The CM procedure is the actual CAS/Mach procedure of the form $\arc_-\arc_{c_1}\arc_{c_2}\arc_+$ while the SA procedure has the simple form $\arc_-\arc_s\arc_+$. Fixing the structure, that is considering the SA or CM procedure, the optimal control problem can be reduced to a simpler optimization problem in finite and small dimension, where the adjoint vector $p$ is not needed.
This is possible in particular because of the parameterization of the control in the feedback form $u_s(x)$, $u_{c_1}(x)$ and $u_{c_2}(x)$ respectively along the arcs $\arc_s$, $\arc_{c_1}$ and $\arc_{c_2}$.
The numerical investigation shows that the SA procedure is better than the actual CM procedure, in terms of time-to-climb and fuel consumption when considering the cost index criterion, but the difference is small. Besides, thanks to the numerical results obtained with the direct collocation method, considering the optimal control problem 
\eqref{eq:ocp_cost_index} with $\phimax = VMO$ and $\psimax=MMO$, we have also that the SA procedure is possibly the best policy.

It is worth to mention that in this article we consider only the {flight path angle} as control variable. This means that the thrust is assumed to be constant. It would be interesting to consider the thrust as a control variable and in this case, we would have a bi-input control system still affine in the control and more complex structures may appear. A preliminary study has been done in \cite{ocam:2017} but the singular perturbation phenomenon has not been taken into account, which makes the numerical resolution much more difficult.


\end{document}